\documentclass[leqno,11pt]{article}
\usepackage{latexsym}
\usepackage{amssymb}
\usepackage{amsfonts}
\usepackage{amsmath}
\usepackage{epsfig}

\makeatletter
\long\def\unmarkedfootnote#1{{\long\def\@makefntext##1{##1}\footnotetext{#1}}}
\makeatother

\newtheorem{definition}{Definition}[section]

\newtheorem{lemma}[definition]{Lemma}

\newtheorem{theorem}[definition]{Theorem}

\newtheorem{proposition}[definition]{Proposition}

\newtheorem{corollary}[definition]{Corollary}

\newtheorem{remark}[definition]{Remark}

\def\o{\Omega}

\def\mo{|\Omega |}
\def\m2{|\Omega | /2}
\def\M2{\frac{|\Omega |}{2}}
\def\u+{u_+^*}

\def\-p{\overline{p}}
\def\W{W^{1,p}(\Omega)}
\def\w0{{W_0^{1,p}(\Omega)}}

\def\R{\mathbb R}
\def\N{\mathbb N}

\def\rn{{{\R}^n}}
\def\T{T_{t}}
\def\Mc{\mathcal{M}}

\newcommand{\hh}{{\cal H}^{n-1}}
\newcommand{\medint}{-\kern  -,395cm\int}
\newcommand{\medintinrigo}{-\kern  -,315cm\int}
\newcommand{\medelle}{-\kern  -,235cm L}
\newcommand{\medellenrigo}{-\kern  -,180cm L}
\newcommand{\qed}{\thinspace\null\nobreak\hfill
\hbox{\vbox{\kern-.2pt\hrule height.2pt
depth.2pt\kern-.2pt\kern-.2pt \hbox to1.8mm {\kern-.2pt\vrule
width.4pt \kern-.2pt\raise1.8mm\vbox to.2pt{} \lower0pt\vtop
to.2pt{}\hfil\kern-.2pt \vrule
width.4pt\kern-.2pt}\kern-.2pt\kern-.2pt \hrule height.2pt
depth.2pt \kern-.2pt}}\par\medbreak}

\title{Well-posed elliptic Neumann problems \\ involving irregular data and domains} 
\numberwithin{equation}{section}

\author{
Angelo Alvino\\
{\it Dipartimento di Matematica e Applicazioni ``R. Caccioppoli"}
\\{\it Universit\`a di Napoli ``Federico II"}
\\ {\it Complesso Monte S.Angelo, Via Cintia, 80126 Napoli, Italy}
\\ {\it e-mail: angelo.alvino@dma.unina.it }
\bigskip
\\
  Andrea Cianchi\\
 {\it Dipartimento di Matematica e Applicazioni per
l'Architettura, Universit\`a di Firenze}\\ {\it Piazza Ghiberti
27, 50122 Firenze, Italy} \\{\it  e-mail: cianchi@unifi.it}
\bigskip
\\
  Vladimir G. Maz'ya \\
 {\it  Department of Mathematical Sciences, M\&O Building}\\ {\it University of Liverpool, Liverpool L69 3BX,
 UK;}\\ and \\
{\it   Department of Mathematics, Link\"oping University, SE-581
83 Link\"oping, Sweden}
\\ {\it e-mail: vlmaz@mai.liu.se}
\bigskip
\\
Anna Mercaldo \\
{\it Dipartimento di Matematica e Applicazioni ``R. Caccioppoli"}
\\{\it Universit\`a di Napoli ``Federico II"}
\\ {\it Complesso Monte S.Angelo, Via Cintia, 80126 Napoli, Italy}
\\{\it e-mail: mercaldo@unina.it }}

\date{}

\begin{document}
\maketitle
%
%
\begin{equation*}
\hbox{{\bf Abstract}} \end{equation*}
 Nonlinear elliptic Neumann problems, possibly in irregular domains
and with data affected by low integrability properties, are taken
into account. Existence, uniqueness and continuous dependence on
the data of generalized solutions are established under a suitable
balance between the  integrability of the datum and the
(ir)regularity of the domain. The latter is described in terms of
isocapacitary inequalities. Applications to various classes of
domains are also presented.
\medskip
\par\noindent
\begin{equation*}
\hbox{{\bf R\'esum\'e}} \end{equation*}
Nous consid\'erons des probl\`emes de Neumann pour des \'equations
elliptiques non lin\'eaires dans domaines \'eventuellement non
r\'eguliers et avec des donn\'ees peu r\'eguli\`eres.
Un \'equilibre entre
l'int\'egrabilit\'e de la donn\'ee et l'(ir)r\'egularit\'e du domaine nous
permet d'obtenir l'existence, l'unicit\'e et la d\'ependance continue
de solutions g\'en\'eralis\'ees. L'irr\'egularit\'e du domaine est d\'ecrite
par des inegalit\'es ``isocapacitaires". Nous donnons aussi des
applications \`a certaines classes de domaines.

\unmarkedfootnote {
}

\section{Introduction and main results}\label{sec1}

The present paper deals with existence, uniqueness and continuous
dependence on the data of solutions  to nonlinear elliptic Neumann
problems having the form
\begin{equation}\label{1}
\begin{cases}
- {\rm div} (a(x,\nabla u)) = f(x)  & {\rm in}\,\,\, \o \\
 a(x,\nabla u)\cdot {\bf  n} =0 &
{\rm on}\,\,\,
\partial \o \,.
\end{cases}
\end{equation}
Here: \par $\o$ is a connected open set in $\rn$, $n\geq 2$,
having finite Lebesgue measure $|\o |$;
\par $a~:~\o~\times~\rn~\rightarrow~ \rn$ is a Carath\'eodory
function;
\par
 $f\in L^q(\o )$ for some $q\in [1, \infty]$ and
satisfies the compatibility condition
\begin{equation}\label{mv}
\int _\o f(x) dx =0.
\end{equation}
\par\noindent
Moreover, $``\cdot "$ stands for inner product in $\rn$, and
${\bf n}$ denotes the outward unit normal on $\partial \o$.
\par
Standard assumptions in the theory of nonlinear elliptic partial
differential equations amount to requiring that there exist an
exponent $p>1$, a function $h\in L^{p'}(\o )$, where
$p'=\frac{p}{p-1}$, and a constant $C$  such that, for a.e. $x\in
\o$ :
\par\noindent
\begin{equation}\label{2}
a(x,\xi)\cdot \xi \geq |\xi |^p \quad \hbox{for} \,\,\, \xi  \in
\rn\,;
\end{equation}
\par\noindent
\begin{equation}\label{3}
|a(x,\xi)| \leq C( |\xi |^{p-1} + h(x)) \quad \hbox{for} \,\,\,
\xi \in \rn\,;
\end{equation}
\par\noindent
\begin{equation}\label{4}
[a(x,\xi)- a(x,\eta)]\cdot (\xi- \eta) > 0 \quad \hbox{for} \,\,\,
\xi , \eta \in \rn\,\,\, \hbox{with} \,\,\,\xi \neq \eta\,.
\end{equation}
\par
The $p$-Laplace equation, corresponding to the choice $a(x,\xi)=
|\xi |^{p-2}\xi$, and, in particular, the (linear) Laplace
equation when $p=2$, can be regarded as prototypal examples on
which our analysis provides new results.
%
\par
 When $\o$ is sufficiently regular, say with a Lipschitz
 boundary, and $q$ is so large that $f$ belongs to the
topological dual of the classical Sobolev space $\W$, namely
$q>\frac{np}{np-n+p}$ if $p<n$, $q>1$ if $p=n$, and $q\geq 1$ if
$p>n$, the existence of a unique (up to additive constants) weak
solution to problem \eqref{1} under \eqref{mv}-\eqref{4} is well
known, and quite easily follows via the Browder-Minthy theory of
monotone operators.
\par
In the present paper, problem \eqref{1} will be set in a more
general framework, where these customary assumptions on $\o$ and
$f$ need not be satisfied. Of course, solutions to \eqref{1} have
to be interpreted in an extended sense in this case. The notion of
solution $u$, called approximable solution throughout this paper,
that will be adopted arises quite naturally in dealing with
problems involving irregular domains and data.
 Loosely speaking, it amounts to demanding
that $u$ be a distributional solution to \eqref{1} which can be
approximated by a sequence of solutions to problems with the same
differential operator and boundary condition, but with regular
right-hand sides. A precise definition can be found in Section
\ref{subsec2.2}. We just anticipate here that an approximable
solution $u$ need not be a Sobolev function in the usual sense;
nevertheless, a generalized meaning to its gradient $\nabla u$ can
be given.
\par Definitions of solutions of
this kind, and other definitions which, a posteriori, turn out to
be equivalent, have been extensively employed in the study of
elliptic Dirichlet problems with a right-hand side $f$ affected by
low integrability properties. Initiated in \cite{Ma0bis, Ma3} and
\cite{St}  in the linear case,
 and in \cite{BG1, BG2} in the nonlinear case, this study
 has been the object of several contributions in the last twenty years, including
\cite{AM, BBGGPV, DallA, DM, DMOP, D, LM, M1, M2, DHM, GIS, FS}. These
investigations have pointed out that, when dealing with
(homogeneous) Dirichlet boundary conditions, existence and
uniqueness of solutions can be established as soon as $f\in
L^1(\o)$, whatever $\o$ is. In fact, the regularity of $\o$ does
not play any role in this case, the underlying reason being that
the level sets of solutions cannot reach $\partial \o$.
 \par The situation is different when Neumann
boundary conditions are prescribed. Actually, inasmuch as the
boundary of the level sets of solutions and $\partial \o$ can
actually overlap, the geometry of the domain $\o$ comes now into
play. We shall prove that problem \eqref{1} is still uniquely
solvable, provided that the (ir)regularity of $\o$ and the
integrability of $f$ are properly balanced. In fact,  even if $f$
highly integrable, in particular essentially bounded, some
regularity on $\o$ has nevertheless to be retained. In the special
case when $\o$ is smooth, or at least with a Lipschitz boundary,
our results overlap with contributions from \cite{AMST, BeGu, Dr,
DV, P,  Pr, R}.
\par Our approach relies upon  isocapacitary inequalities, which have recently been shown in \cite{CM1}  to provide
suitable information on the regularity of the domain $\o$ in the
study of problems of the form \eqref{1}.
   In fact, isocapacitary inequalities turn out to be more effective than the more popular
isoperimetric inequalities in this kind of applications. The use
of the standard isoperimetric inequality in the study of
 elliptic Dirichlet problems, and of  relative isoperimetric
 inequalities in the study of
  Neumann problems, was introduced in
 \cite{Ma0bis, Ma3}. The  isoperimetric inequality
 was also independently employed in \cite{Ta1, Ta3} in the
 proof of symmetrization principles for solutions to Dirichlet
 problems. Ideas from these papers have been developed in a rich literature, including
 \cite{Al, AFLT, ALT, Kesavan}. Specific contributions to the study of Neumann
 problems are  \cite{AMT, Be, Celliptic,
 FeEdy, Fe, MS1, MS2}. We refer to \cite{Ke, Tr, Va}
 for an exhaustive bibliography on these topics.
%
%
%
\par
 The relative
isoperimetric inequality in $\o$
tells us that
\begin{equation}\label{6}
\lambda (|E|) \leq P(E; \o) \qquad \quad \hbox{for every
measurable set $E \subset \o$ with $|E| \leq \m2$},
\end{equation}
where $P(E; \o)$ denotes the perimeter of a measurable set $E$
relative to $\o$, and  $\lambda : [0 , \m2 ] \rightarrow [0,
\infty )$ is the isoperimetric function of $\o$.
\par\noindent
Replacing the relative perimeter by a suitable $p$-capacity on the
right-hand side of \eqref{6} leads to the isocapacitary inequality
in $\o$. Such inequality reads
\begin{equation}\label{7}
\nu _p (|E|) \leq C_p(E,G) \quad \hbox{for every measurable sets
$E\subset G \subset \o$ with $|G| \leq |\o |/2$},
\end{equation}
where $C_p(E,G)$ is the $p$-capacity of the condenser $(E;G)$
relative to $\o$, and ${\nu _p}: [0 , \m2 ] \rightarrow [0, \infty
]$ is the isocapacitary function of $\o$.
\par Precise
definitions concerning perimeter and capacity, together with their
properties entering in our discussion, are given in Section
\ref{subsec2.3}. Let us emphasize that although \eqref{6} and
\eqref{7} are essentially equivalent for sufficiently smooth
domains $\o$, the isocapacitary inequality \eqref{7} offers, in
general, a finer description of the regularity of bad domains
$\o$. Accordingly, our main results will be formulated and proved
in terms of the function $\nu _p$. Their counterparts involving
$\lambda$ will be derived as corollaries  - see Section
\ref{sec5}. Special instances of bad domains and data will
demonstrate that the use of $\nu _p$ instead of $\lambda$ can
actually lead to stronger conclusions in connection with
existence, uniqueness and continuous dependence on the data of
solutions to problem \eqref{1}.
\par
Roughly speaking, the faster   the function $\nu _p(s)$ decays to
$0$ as $s \to 0^+$, the worse is the domain $\Omega$, and,
obviously, the smaller  is $q$, the worse is $f$. Accordingly, the
spirit of our results is that problem \eqref{1} is actually
well-posed, provided that $\nu _p(s)$  does not decay to $0$ too
fast as $s \to 0^+$, depending on how small $q$ is. Our first
theorem provides us with conditions for the unique solvability (up
to additive constants) of \eqref{1} under the basic assumptions
\eqref{mv}--\eqref{4}.
\begin{theorem}\label{teo4bis}
Let $\o$ be an open connected subset of $\rn$, $n \geq 2$, having
finite measure. Assume that  $f\in L^q(\o )$  for some $q\in [1,
\infty ]$ and satisfies \eqref{mv}.
Assume that \eqref{2}-\eqref{4} are fulfilled, and that either
\par\noindent
{\rm (i)} $1<q \leq \infty $ and
\begin{equation}\label{401}
\int _0^{\m2}\bigg(\frac{s}{\nu _p(s)}\bigg)^{\frac {q'}{p}} ds <
\infty\,,
\end{equation}
or
\par\noindent
{\rm (ii)} $q=1$ and
\begin{equation}\label{402}
\int _0^{\m2}\bigg(\frac{s}{\nu _p(s)}\bigg)^{\frac{1}{p}}
\frac{ds}{s} < \infty\,.
\end{equation}
%
%
%
Then there exists a unique (up to additive constants) approximable
solution to problem \eqref{1}.
\end{theorem}
\par The second main result of this paper is concerned with the
case when the differential operator in \eqref{1} is not merely
strictly monotone in the sense of \eqref{4}, but fulfils the
strong monotonicity assumption that, for a.e. $x\in \o$,
\begin{equation}\label{5}
[a(x,\xi)- a(x,\eta)]\cdot (\xi- \eta) \geq
\begin{cases}
C |\xi - \eta|^p  & {\rm if}\,\,\, p\geq 2 \\
C \displaystyle {\frac{|\xi - \eta|^2}{(|\xi| + |\eta |)^{2-p}} }&
{\rm if}\,\,\, 1<p<2 \,,
\end{cases}
\end{equation}
for some positive constant $C$ and for $\xi, \eta \in \rn$.  In
addition to the result of Theorem \ref{teo4bis}, the continuous
dependence of the solution to \eqref{1} with respect to $f$ can be
established under the reinforcement of \eqref{4} given by
\eqref{5}. In fact,  when \eqref{5}  is in force, a partially
different approach can be employed, which also simplifies the
proof of the statement of Theorem \ref{teo4bis}.
%
\par Observe that, in particular, assumption
\eqref{5} certainly holds provided that, for a.e. $x\in \o$, the
function $a(x, \xi )= (a_1(x, \xi ), \dots , a_n(x, \xi ))$ is
differentiable with respect to $\xi$, vanishes for $\xi =0$, and
satisfies the ellipticity condition \begin{equation*} \sum
_{i,j=1}^{n} \frac{\partial a_i}{\partial \xi _j}(x, \xi ) \eta _i
\eta _j \geq C |\xi |^{p-2}|\eta |^2 \qquad \hbox{for
 $\xi, \eta \in \rn$},
 \end{equation*}
 for some positive constant $C$.
\begin{theorem}\label{teo4}
Let $\o$, $p$, $q$ and $f$ be as in Theorem \ref{teo4bis}. Assume
that \eqref{2}, \eqref{3}  and  \eqref{5} are fulfilled. Assume
that either $1<q \leq \infty $ and \eqref{401} holds, or $q=1$ and
\eqref{402} holds.
%
%
\par\noindent
Then there exists a unique (up to additive constants) approximable
solution to problem \eqref{1}  depending continuously on the
right-hand side of the equation. Precisely, if $g$ is another
function from $L^q(\o )$ such that $\int _\o g(x) dx =0$, and $v$
is the solution to \eqref{1} with $f$ replaced by $g$, then
\begin{equation}\label{403}
\|\nabla u - \nabla v\|_{L^{p-1}(\o)} \leq C \|f -
g\|_{L^{q}(\o)}^{\frac 1r} \big(\|f\|_{L^{q}(\o)}+
\|g\|_{L^{q}(\o)}\big)^{\frac 1{p-1} - \frac 1r}
\end{equation}
for some constant $C$ depending on $p$, $q$ and on the left-hand
side  either of   \eqref{401} or  \eqref{402}. Here, $r = \max
\{p, 2\}$.
\end{theorem}
\par
Let us notice that the balance condition between $q$ and $\nu _p$
in Theorems \ref{teo4bis} and \ref{teo4} requires a separate
formulation according to whether $q>1$ or $q=1$. In fact,
assumption \eqref{402}  is a qualified version  of the limit as
$q\to 1^+$ of \eqref{401}. This is as a consequence of the
different a priori (and continuous dependence) estimates upon
which Theorems \ref{teo4bis} and \ref{teo4} rely. Actually,
$L^1(\o)$ is a borderline space, and when $f \in L^1(\o )$  the
natural sharp estimate  involves a weak type (i.e. Marcinkiewicz)
norm of the gradient of the solution $u$. Instead, when  $f \in
L^q(\o)$ with $q>1$, a strong type   (i.e. Lebesgue) norm comes
into play in a sharp bound for the gradient of $u$. This gap is
intrinsic in the problem, as witnessed by the basic case of the
Laplace (or p-Laplace) operator in a smooth domain.
\par
\smallskip
 The paper is organized as follows. In Section \ref{sec2}
we collect definitions and basic properties concerning functions
spaces of measurable (Subsection \ref{subsec2.1}) and weakly
differentiable functions (Subsection \ref{subsec2.1'}), solutions
to problem \eqref{1} (Subsection \ref{subsec2.2}), perimeter and
capacity (Subsection \ref{subsec2.3}). Section \ref{sec3} is
devoted to the proof of Theorem \ref{teo4bis}, which is
accomplished in Subsection \ref{subsec3.2}, after deriving the
necessary a priori estimates in Subsection \ref{subsec3.1}.
Continuous dependence estimates under the  strong monotonicity
assumption \eqref{5} are established in Subsection \ref{subsec4.1}
of Section \ref{sec4}; they are a key step in the proof of Theorem
\ref{teo4} given in Subsection \ref{subsec4.2}. Finally, Section
\ref{sec5} contains applications of our results to special domains
and classes of domains. Versions of Theorems \ref{teo4bis} and
\ref{teo4} involving the isoperimetric function are also
preliminarily stated. With their help, the advantage of the use of
isocapacitary inequalities instead of isoperimetric inequalities
is demonstrated in concrete examples.

\section{Background and preliminaries}\label{sec2}

\subsection{Rearrangements and rearrangement invariant spaces}\label{subsec2.1}

Let us denote by $\Mc (\Omega )$  the set of measurable functions
in $\Omega$, and let $u\in \Mc (\Omega )$. The distribution
function $\mu_u: [0, \infty)\rightarrow [0, \infty)$ of $u$ is
defined as
\begin{equation}\label{mu}
\mu_u (t) = |\{x\in \o:÷|u(x)|\ge t\}|, \qquad \hbox{for $t\ge
0$}.
\end{equation}
The decreasing rearrangement $u^*: [0, |\o|]\rightarrow [0,
\infty]$ of $u$ is given by
\begin{equation}\label{decrear}
u^*(s)=\sup\,\{t\ge 0 :÷\mu_u(t)\ge s\}|,  \qquad \hbox{for }÷s\in
[0, |\o|].
\end{equation}
We also define $u_*: [0, |\o|]\rightarrow [0, \infty ]$, the
increasing rearrangement of $u$,  as
$$
u_*(s)=u^*(|\o|-s), \qquad \hbox{for }÷s\in [0, |\o|].
$$
The operation of decreasing rearrangement is neither additive nor
subadditive. However,
\begin{equation}\label{somma}
(u+v)^*(s)\le u^*(s/2)+v^*(s/2), \qquad \hbox{for }÷s\in [0,
|\o|],
\end{equation}
for any $u, v \in \Mc (\o )$, and hence, via Young's inequality,
\begin{equation}\label{prodotto}
(uv)^*(s)\le u^*(s/2)v^*(s/2), \qquad \hbox{for }÷s\in [0, |\o|].
\end{equation}
\par\noindent
A basic property of  rearrangements is the Hardy-Littlewood
inequality, which tells us that
\begin{equation}\label{HL}
\int_0^{|\o|}u^*(s)v_*(s)\, ds\le \int_{\o}|u(x)v(x)|\, dx \le \int_0^{|\o|}u^*(s)v^*(s)\, ds\
\end{equation}
for any  $u, v \in \Mc (\o )$.
\par
A rearrangement invariant (r.i., for short) space $X(\o)$ on $\o$
is a Banach function space,  in the sense of Luxemburg, equipped
with a norm $\|\cdot \|_{X(\o)}$ such that
\begin{equation}\label{norm}
\|u\|_{X(\o)}=\|v\|_{X(\o)} \qquad \quad \hbox{whenever
$u^*=v^*$}.
\end{equation}
Since we are assuming that $|\o |<\infty$, any r.i. space $X(\o )$
fulfills
$$L^\infty (\o ) \to X(\o ) \to L^1(\o )\,,$$
where the arrow ``$\to$"  stands for continuous embedding.
\par
Given any r.i. space $X(\o )$, there exists a unique r.i. space
$\overline X(0, |\o|)$, the representation space of $X(\o )$ on
$(0, |\o|)$, such that
\begin{equation}\label{2.8}
\|u\|_{X(\o)}=\|u^*\|_{\overline X(0, |\o|)}
\end{equation}
for every $u \in X(\o)$. A characterization of the norm $\|\cdot
\|_{\overline X(0, |\o|)}$ is available (see \cite[Chapter 2,
Theorem 4.10 and subsequent remarks]{BS}). However, in our
applications, an expression for $\overline X(0, |\o|)$ will be
immediately derived via basic properties of rearrangements. In
fact, besides the standard Lebesgue spaces, we shall only be
concerned with Lorentz and Marcinkiewicz type spaces. Recall that,
given $\sigma, \varrho \in (0, \infty)$, the Lorentz space
$L^{\sigma, \varrho}(\o)$ is the set of all functions $u \in \Mc
(\o )$ such that the quantity
\begin{equation}\label{normlor}
\|u\|_{L^{\sigma, \varrho}(\o)}=\left (\int_0^{|\o|}
(s^\frac{1}{\sigma}u^*(s))^\varrho\frac{ds}{s}\right )^{1/\varrho}
\end{equation}
is finite. The expression $\|\cdot \|_{L^{\sigma, \varrho}(\o)}$
is an (r.i.) norm if and only if $1 \leq \varrho \le \sigma $.
When $\sigma \in (1, \infty )$ and $\varrho \in [1, \infty )$, it
is always equivalent to the norm obtained on replacing $u^*(s)$ by
$\frac{1}{s}\int_0^su^*(r) dr$ on the right-hand side of
\eqref{normlor}; the space $L^{\sigma, \varrho}(\o)$, endowed with
the resulting norm, is an r.i. space. Note that $L^{\sigma,
\sigma}(\o)=L^{\sigma}(\o) $ for $\sigma
>0$. Moreover, $L^{\sigma, \varrho _1}(\o) \subsetneqq L^{\sigma, \varrho
_2}(\o)$ if $\varrho _1 < \varrho _2$, and, since $|\Omega
|<\infty$, $L^{\sigma _1, \varrho _1}(\o) \subsetneqq L^{\sigma
_2, \varrho _2}(\o)$ if $\sigma _1 > \sigma _2$ and $\varrho _1 ,
\varrho _2 \in (0,  \infty].$
\par
Let $\omega ÷: (0,|\o|)\rightarrow (0, \infty)$ be a bounded
non-decreasing function. The Marcinkiewicz space $M_{\omega }(\o)$
associated with $\omega$ is the set of all functions $u\in \Mc (\o
)$ such that the quantity
\begin{equation}\label{normmarc}
\|u\|_{M_{\omega}(\o)}=\sup_{s\in (0, |\o|)} \omega (s) u^*(s)
\end{equation}
is finite. The expression \eqref{normmarc} is equivalent to a
norm,  which makes $M_{\omega}(\o)$  an r.i. space, if and only if
$\sup _{s\in (0, |\Omega |)}\frac {\omega(s)}s \int
_0^s\frac{dr}{\omega (r)}<\infty$\,.
\subsection{Spaces of Sobolev type}\label{subsec2.1'}

Given any $p\in [1, \infty]$, we denote by
$W^{1,p}(\o)$ the standard Sobolev space, namely
$$
W^{1,p}(\o)=\{ u\in L^p(\o): \,\, \hbox{$u$ is weakly
differentiable in $\o$ and $|\nabla u|\in L^p(\o)$} \}.
$$
The space $W_{loc}^{1,p}(\o)$ is defined analogously, on replacing
$L^p(\o)$ by  $L_{loc}^p(\o)$ on the right-hand side.
\par
Given any $t>0$,  let  $T_{t} : \R \rightarrow \R$ be the function
given by
\begin{equation}\label{502}
T_{t} (s) = \begin{cases}
s & {\rm if}\,\,\, |s|\leq t \\
 t \,{\rm sign}(s) &
{\rm if}\,\,\, |s|>t \,.
\end{cases}
\end{equation}
For $p \in [1, \infty ]$, we set
\begin{equation}\label{503}
W^{1,p}_T(\o) = \left\{u: \hbox{$u\in \Mc (\o )$ and  $T_{t}(u)
\in W^{1,p}(\o)$ for every $t >0$} \right\}\,.
\end{equation}
The space  $W_{T, loc}^{1,p}(\o)$ is defined accordingly, on
replacing $W^{1,p}(\o)$  by $W_{loc}^{1,p}(\o)$ on the right-hand
side of \eqref{503}. If $u\in W_{T, loc}^{1,p}(\o)$, there exists
a (unique) measurable function  $Z_u : \o \to \rn$ such that
\begin{equation}\label{504}
\nabla \big(T_{t}(u)\big) = \chi _{\{|u|<k\}}Z_u  \qquad \quad
\hbox{a.e. in $\o$}
\end{equation}
for every $t>0$ \cite[Lemma 2.1]{BBGGPV}. Here $\chi _E$ denotes
the characteristic function of the set $E$. One has that $u \in
W_{loc}^{1,p}(\o)$ if and only if $u \in W_{T, loc}^{1,p}(\o)\cap
L_{loc}^p(\o )$ and  $Z_u \in L_{loc}^p(\o , \rn )$, and, in this
case, $Z_u= \nabla u$. An analogous property holds provided that
``loc" is  dropped everywhere. In what follows, with abuse of
notation, for every $u\in W_{T, loc}^{1,p}(\o)$  we denote $Z_u$
by $\nabla u$.
\par
Given $p\in (0, \infty]$, define
$$
V^{1,p}(\o)= \left\{u: u\in W_{T, loc}^{1,1}(\o) \hbox{ and } |\nabla u|\in L^p(\o)  \right\}\,.
$$
Note that, if  $p \geq 1$, then
$$
V^{1,p}(\o)= \left\{u: u\in W_{ loc}^{1,1}(\o) \hbox{ and } |\nabla u|\in L^p(\o)  \right\}\,,
$$
a customary space of weakly differentiable functions. Moreover, if
 $p\ge 1$,  the set $\o$ is connected, and $B$ is any ball such that  $\overline{B}\subset \o$,
 then $V^{1,p}(\o)$ is a Banach space equipped with the norm
$$
\|u\|_{V^{1,p}(\o)}=\|u\|_{L^{p}(B)}+\|\nabla u\|_{L^{p}(\o)}.
$$
Note that, replacing $B$ by another ball results  an equivalent
norm. The topological dual of $V^{1,p}(\o)$ will be denoted by
$(V^{1,p}(\o))'$.

Given any ball $B$ as above, define the subspace $V^{1,p}_B(\o)$
of $V^{1,p}(\o)$ as
$$
V^{1,p}_B(\o)= \left\{u\in V^{1,p}(\o): \,\, \int_B u\, dx=0
\right\}\,.
$$
\begin{proposition}\label{(B)}
Let $p\in [1, \infty]$. Let $\o$ be a connected open set in $\rn$
having finite measure, and let $B$ be any ball such that
$\overline{B}
 \subset \o$. Then the quantity
\begin{equation}\label{propB}
\|u\|_{V_B^{1,p}(\o)}=\|\nabla u\|_{L^{p}(\o)}
\end{equation}
defines a norm in $V_B^{1,p}(\o)$ equivalent to $\| \cdot
\|_{V^{1,p}(\o)}$. Moreover, if $p\in (1, \infty )$, then
$V_B^{1,p}(\o)$, equipped with this norm, is a separable and
reflexive Banach space.
\end{proposition}
\noindent {\bf Proof, sketched}. The only nontrivial property that
has to be checked in order to show that $\|\cdot
\|_{V_B^{1,p}(\o)}$ is actually a norm is the fact that
$\|u\|_{V_B^{1,p}(\o)} = 0$ only if $u=0$. This is a consequence
of the Poincare type inequality which tells us that, for every
smooth open set $\o '$ such that  $\overline{B}\subset \o '$ and
$\overline{\o '} \subset  \o $,
\begin{equation}\label{propB'}
\|u\|_{L^{p}(\o ')}\le C \|\nabla u\|_{L^{p}(\o ')}\,
\end{equation}
for some constant $C=C(p, \o ', |B|)$ and for every $u\in
V_B^{1,p}(\o)$ (see e.g. \cite[Chapter 4]{Z}).  The same
inequality plays a role in showing that $V_B^{1,p}(\o)$, equipped
with the norm $\|\cdot \|_{V_B^{1,p}(\o)}$, is complete. When
$p\in (1, \infty )$, the separability and the reflexivity of
$V_B^{1,p}(\o)$ follow via the same argument as for the standard
Sobolev space $W^{1,p}(\o)$, on making use of the fact that the
map $L: V_B^{1,p}(\o)\rightarrow (L^p(\o))^n$ given by $Lu=\nabla
u$ is an isometry of $V_B^{1,p}(\o)$ into $(L^p(\o))^n$, and that
$ (L^p(\o))^n$ is a separable and reflexive Banach space.
\medskip



\subsection{Solutions}\label{subsec2.2}

When $f\in (V^{1,p}(\o ))'$, and \eqref{mv}--\eqref{3} are in
force, a standard notion of solution to problem \eqref{1} is that
of weak solution. Recall that a function $u\in V^{1, p}(\o)$ is
called a weak solution to \eqref{1} if
\begin{equation}\label{weaksol}
\int_\o a(x, \nabla u)\cdot \nabla \Phi \, dx=\int_\o f\Phi\, dx\,
\qquad \quad \hbox{for every $\Phi \in V^{1, p}(\o)$.}
\end{equation}
\par
\noindent An  application of the Browder-Minthy theory for
monotone operators, resting upon Proposition \ref{(B)}, yields the
following  existence and uniqueness result. The proof can be
accomplished along the same lines as in \cite[Porposition 26.12
and Corollary 26.13]{Zeidler}. We omit the details for brevity.
\begin{proposition}\label{weak}
Let $p\in (1, \infty )$ and let $\o$ be a bounded connected open
set in $\rn$ having finite measure. If $f\in (V^{1,p}(\o ))'$,
then under assumptions \eqref{mv}--\eqref{4} there exists a unique
(up to additive constants) weak solution $u\in V^{1, p}(\o)$ to
problem \eqref{1}.
\end{proposition}

The definition of weak solution does not fit the case when
$f\notin (V^{1,p}(\o ))'$, since the right-hand side of
\eqref{weaksol} need not be well-defined. This difficulty can be
circumvented on restricting the class of test functions $\Phi$ to
$W^{1, \infty}(\o)$, for instance. This leads to a counterpart, in
the Neumann problem setting, of the classical definition of
solution to the Dirichlet problem in the sense of distributions.
It is however well-known \cite{Se} that such a  class of test
functions  may be too poor for the solution to be uniquely
determined, even under an appropriate monotonicity assumption as
\eqref{4}.
%

In order to overcome this drawback, we adopt a definition of
solution, in the spirit e.g. of \cite{DallA} and \cite{DM},
obtained in the limit from solutions to approximating problems
with regular right-hand sides. The idea behind such a definition
is that the additional requirement of being approximated by
solutions to regular problems identifies a distinguished proper
distributional solution to problem \eqref{1}.
Specifically, if $\o $ is an open set in $\rn$ having finite
measure,
and $f\in L^q(\o)$ for some $q\in [1, \infty ]$ and fulfills
\eqref{mv}, then a function $u\in V^{1,p-1}(\o)$ will be called an
\emph{approximable solution} to problem \eqref{1} under
assumptions \eqref{2} and \eqref{3} if:
\par\noindent
(i)
\begin{equation}\label{231}
\int_\o a(x, \nabla u)\cdot \nabla \Phi \, dx=\int_\o f\Phi\, dx\,
\qquad \quad \hbox{for every $\Phi \in W^{1, \infty}(\o)$,}
\end{equation}
and
\par\noindent
(ii) a sequence $\{f_k\}\subset  L^q(\o) \cap (V^{1,p}(\o ))'$
exists such that $\int _\o f_k(x)dx=0$ for $k \in \N$,
$$f_k
\rightarrow f \qquad \hbox{in $L^q(\o)$},$$
and the sequence of
weak solutions $\{u_k\}\subset V^{1, p}(\o)$ to problem \eqref{1},
with $f$ replaced by $f_k$, satisfies $$u_k\rightarrow u \quad
\hbox{ a.e. in $\o$}.$$

A few brief comments about this definition  are in order.
Customary counterparts of such a definition for Dirichlet problems
\cite{DallA, DM} just amount to (a suitable version of) property
(ii). Actually, the existence of a generalized gradient of the
limit function  $u$, in the sense of \eqref{504}, and the fact
that $u$ is a distributional solution
 directly follow from analogous properties of the
approximating solutions $u_k$. This is due to the fact that,
whenever $f\in L^1(\Omega )$, a priori estimates in suitable
Lebesgue spaces for the gradient of approximating solutions to
homogeneous Dirichlet problems
are available, irrespective of whether $\Omega$ is regular or not.
As a consequence, one can pass to the limit in the equations
fulfilled by $u_k$, and hence  infer that $u$ is a distributional
solution to the original Dirichlet problem.
When Neumann
problems are taken into account, the existence of a generalized
gradient of $u$ and the validity of (i) is not guaranteed anymore,
inasmuch as a priori estimates for $|\nabla u_k|$ depend on the
regularity of $\o$. The membership of $u$ in $V^{1, p-1}(\o)$ and
equation (i) have consequently to be included as part of the
definition of solution.
\par
Let us also mention that the definition of approximable solution
can be shown to be equivalent to other definitions patterned on
those of entropy solution \cite{BBGGPV} and of renormalized
solution \cite{LM} given for Dirichlet problems.
\subsection{Perimeter and capacity}\label{subsec2.3}

The isoperimetric function $\lambda : [0, \m2 ] \to [0, \infty )$
of $\o$
 is defined
as
%
\begin{equation}\label{-512}
\lambda (s) = \inf \{P(E, \o ):  s \leq |E| \leq \m2 \} \qquad
\quad \hbox{for $s\in [0, \m2 ]$}\,.
\end{equation}
Here, $P(E; \o)$ is the perimeter of $E$ relative to $\o$, which
agrees with $\hh (\partial ^M E \cap \o)$, where $\hh$ denotes the
$(n-1)$-dimensional Hausdorff measure, and $\partial ^M E$ stands
for the essential boundary of $E$ (see e.g. \cite{AFP, Ma2}).
\par
The relative isoperimetric inequality \eqref{6}
 is a straightforward consequence of definition \eqref{-512}. On the other hand,
 the isoperimetric function $\lambda$ is  known only for
very special domains, such as balls \cite{Ma2, BuZa} and convex
cones \cite{LP}. However, various qualitative and quantitative
properties of $\lambda$ have been investigated, in view of
applications to Sobolev inequalities \cite{HK, Ma0, Ma2, MP},
eigenvalue estimates \cite{Cheeger, Celliptic, Gallot}, a priori
bounds for solutions to Neumann problems (see the references in
Section \ref{sec1}).
\par\noindent
In particular, the function $\lambda$ is  known to be strictly
positive in $(0, \m2 ]$ when $\o$ is connected \cite[Lemma
3.2.4]{Ma2}. Moreover, the asymptotic behavior of $\lambda (s)$ as
$s \to 0^+$ depends on the regularity of the boundary of $\o$. For
instance, if $\o$ has a Lipschitz boundary, then
\begin{equation}\label{-515}
\lambda (s) \approx s^{1/n'} \qquad \quad \hbox{as $s \to 0^+$}
\end{equation}
\cite[Corollary 3.2.1/3]{Ma2}. Here, and in what follows, the
relation $\approx$ between two quantities
 means that the relevant quantities are bounded by each
other up to multiplicative constants. The asymptotic behavior of
the function $\lambda$  for sets having an
 H\"older continuous boundary in the plane was established  in \cite{Crelative}. More general
 results
  for sets in $\rn$ whose boundary has
 an arbitrary modulus of continuity follow from \cite{La}.
 Finer asymptotic estimates for $\lambda$ can be derived under
 additional assumptions on $\partial \o$ (see e.g. \cite{CY, Cmoser}).
\par\noindent
\par
\medskip
The approach of the present paper relies upon estimates for the
Lebesgue measure of subsets of $\o$ via their relative
condenser capacity instead  of their relative perimeter. 
Recall that the standard $p$-capacity of a set $E\subset \o $ can
be defined for $p\geq 1$ as
\begin{equation}\label{-506}
C_p(E) = \inf \left\{ \int _\o |\nabla u|^p\,dx : u\in
W^{1,p}_0(\o), u \geq 1 \,\,\, \hbox{in some neighbourhood
of}\,\,\, E\right\},
\end{equation}
where $W^{1,p}_0(\o)$ denotes the closure in $W^{1,p}(\o)$
of the set of smooth compactly supported functions in $\o $. A
property concerning the pointwise behavior of functions is said to
hold $C_p$-quasi everywhere in $\o$, $C_p$-q.e. for short,  if
it is fulfilled outside a set of $p$-capacity zero.
\par\noindent
Each function $u \in W^{1,p}(\o)$ has a representative
$\widetilde{u}$,  called the precise representative, which is
$C_p$-quasi continuous, in the sense that for every $\varepsilon
>0$, there exists a set $A \subset \o$, with $C_p(A) <
\varepsilon$, such that $f_{|\o \setminus A}$ is continuous in
$\o \setminus A$. The function $\widetilde{u}$ is unique, up
to subsets of $p$-capacity zero. In what follows, we assume that
any function $u \in W^{1,p}(\o)$ agrees with its precise
representative.
\par
A standard result in the theory of capacity tells us that, for
every set $E\subset \o$,
\begin{equation}\label{-507}
C_p(E) = \inf \left\{ \int _\o |\nabla u|^p\,dx : u\in
W^{1,p}_0(\o), u \geq 1 \,\,\, \hbox{$C_p$-q.e. in}\,\,\,
E\right\}
\end{equation}
-- see e.g. \cite[Proposition 12.4]{Da} or \cite[Corollary
2.25]{MZ}. In the light of  \eqref{-507}, we adopt the following
definition of capacity of a condenser. Given sets $E \subset G
\subset \o$, the capacity $C_p(E,G)$  of the condenser $(E, G)$
relative to $\o$ is defined as
\begin{equation}\label{-508}
C_p(E,G) = \inf \left\{ \int _\o |\nabla u|^p\,dx : u\in
W^{1,p}(\o), \hbox{$u \geq 1$ $C_p$-q.e. in $E$ and $u \leq 0$
$C_p$-q.e. in $\o \setminus G$} \right\}\,.
\end{equation}
\par\noindent
Accordingly, the $p$-isocapacitary function ${\nu _p}: [0, |\o
|/2) \to [0, \infty )$ of $\o $ is  given by
\begin{multline}\label{-509}
\nu _p(s) = \inf \left\{C_p(E,G): \hbox{$E$ and $G$ are
measurable subsets of $\o$ such that}\right.\\
\left. \hbox{$E \subset G \subset \o$, $s \leq |E|$ and
$|G|\leq |\o |/2$} \right\}\qquad \hbox{for $s\in [0, |\o
| /2)$}.
\end{multline}
The function $\nu _p$ is clearly non-decreasing.  In what follows, we
shall always deal with the left-continuous representative of
$\nu _p$, which, owing to the monotonicity of $\nu _p$, is pointwise
dominated by the right-hand side of \eqref{-509}.
\par
The isocapacitary inequality  \eqref{7} immediately follows from
definition \eqref{-509}. The point is again to get information
about the behavior of $\nu _p (s)$ as $s \to 0^+$. Such a behavior
is known to related, for instance, to validity of Sobolev
embeddings for $V^{1,p}(\o)$ -- see  \cite{Ma2, MP}, where further
results concerning $\nu _p$ can also be found. In particular, a
slight variant of the results of \cite[Section 8.5]{MP} tells us
that
\begin{equation}\label{2001}
V^{1,p}(\Omega ) \to L^\sigma (\Omega )
\end{equation}
  if and
only if either  $1 \leq p \leq \sigma < \infty$ and
\begin{equation}\label{2002} \sup _{0 < s <\m2} \frac
{s^{\frac p\sigma}}{\nu _p (s)} < \infty\,,
\end{equation}
or $1 \leq \sigma \leq p$ and
\begin{equation}\label{2002'}
\int _0^{\m2} \bigg(\frac{s^{p/\sigma}}{\nu _p (s)}\bigg)^{\frac
\sigma{p-\sigma}}\frac {ds}{s} < \infty\,.
\end{equation}

%
%
\par
As far as relations between $\lambda$ and $\nu _p$ are concerned,
given any connected open set $\Omega$ with finite measure one has
that
\begin{equation}\label{2003}
\nu _1 (s) \approx \lambda (s), \quad \hbox{as $s \to 0^+$,}
\end{equation}
as shown by an easy variant of \cite[Lemma 2.2.5]{Ma2}. When
$p>1$, the functions $\lambda$ and $\nu _p$ are related by
\begin{equation}\label{-514}
\nu _p(s) \geq \bigg(\int _s^{\m2} \frac{dr}{\lambda
(r)^{p'}}\bigg)^{1-p}, \qquad \quad \hbox{for $s\in (0, \m2 )$}
\end{equation}
\cite[Proposition 4.3.4/1]{Ma2}. Hence, in particular, $\nu _p$ is
strictly positive in $(0, \m2 )$ for every connected open set
having finite measure, and
\begin{equation}\label{-509bis}
\lim _{s\to \m2 ^-} \nu _p(s)= \infty \,.
\end{equation}
\par
 A
reverse inequality in \eqref{-514}  does not hold in general, even
up to a multiplicative constant. This accounts for the fact that
the results on problem \eqref{1} which can be derived in terms of
$\nu _p$ are stronger, in general, than those resting upon
$\lambda$. However, the two sides of \eqref{-514} are equivalent
when $\o $ is sufficiently regular. This is the case, for
instance, if $\o$ is bounded and has a Lipschitz boundary. In this
case, combining  \eqref{-515} and \eqref{-514}, and choosing small
concentric balls as sets $E$ and $G$ to estimate the right-hand
side in definition \eqref{-509} easily show that
\begin{equation}\label{-516}
\nu _p(s) \approx s^{\frac{n-p}{n}} \qquad \quad \hbox{as $s \to
0^+$}\,,
 \end{equation}
 if  $p\in [1,n)$, whereas
\begin{equation}\label{-516bis}
\nu _n(s) \approx \Big( \log \frac{1}{s}\Big)^{1-n} \qquad \quad
\hbox{as $s \to 0^+$}\,.
 \end{equation}



\section{Strictly monotone operators}\label{sec3}

\subsection{A priori estimates}\label{subsec3.1}
In view of their use in the proofs of Theorems \ref{teo4bis} and
\ref{teo4}, we collect here a priori estimates for the solution
$u$ to problem \eqref{1} and for its gradient $\nabla u$, under
assumptions \eqref{2}--\eqref{4}. Both pointwise estimates for
their decreasing rearrangements, and norm estimates are presented.
Our results are stated for weak solutions to \eqref{1} under the
assumption that $f\in (V^{1,p}(\o ))'$, this being sufficient for
them to be applied to the approximating problems. We emphasize,
however, that  these results continue to hold for approximable
solutions when $f\notin (V^{1,p}(\o ))'$, as it is easily shown on
 adapting the  approximation arguments that will be
exploited in the proof of Theorem \ref{teo4bis}. Thus, the results
of the present section can also be regarded as  regularity results
 for approximable solutions to problem \eqref{1}.
\par
We begin with  estimates for $u$, which are contained in Theorems
\ref{teoA} and \ref{teoC} below. In what follows, we set
\begin{equation}\label{-520bis}
{\rm med } (u) = \sup \{t\in \R : |\{u>t\}| \geq \m2\}\,,
\end{equation}
the median of $u$. Hence, if
\begin{equation}\label{-520}
{\rm med } (u) = 0\,,
\end{equation}
then
\begin{equation}\label{-521}
|\{u>0 \}| \leq \m2 \qquad \hbox{and} \qquad |\{u <0\}| \leq
\m2\,.
\end{equation}
Moreover, we adopt the notation $u_+ = \frac{|u|+u}{2}$ and
$u_-=\frac{|u|-u}{2}$ for the positive and the negative part of a
function $u$, respectively.
\begin{theorem}\label{teoA}
Let $\o$, $p$ and $a$ be as in Theorem \ref{teo4bis}.
Assume that $f\in L^1(\o )\cap (V^{1,p}(\o ))'$ and fulfills
\eqref{mv}. Let $u$ be the weak solution to problem \eqref{1} such
that ${\rm med } (u) =0$. Then
\begin{equation}\label{A}
u_\pm ^* (s) \leq \int _s^{\m2} \bigg(\int _0^r f_\pm ^* (\rho )
d\rho\bigg)^{\frac{1}{p-1}}\,d(-D\nu _p^{\frac{1}{1-p}})(r), \qquad
\hbox{for  $s\in (0, \m2 )$.}
\end{equation}
Here, $D\nu _p^{\frac{1}{1-p}}$ denotes the derivative
in the sense of measures of the non-increasing function
$\nu_p^{\frac{1}{1-p}}$.
\end{theorem}

\begin{theorem}\label{teoC}
Let $\o$, $p$ and $a$ be as in Theorem \ref{teo4bis}. Assume that
 $f\in L^q(\o)\cap (V^{1,p}(\o ))'$ for some $q\in [1, \infty]$ and fulfills \eqref{mv}. Let
$u$ be the weak solution to problem \eqref{1} such that ${\rm med
} (u) =0$. Let $\sigma \in (0, \infty )$. Then, there exists a
constant $C$ such that
\begin{equation}\label{Co}
\|u \|_{L^\sigma (\o )} \leq C \|f\|_{L^q(\o )}^{\frac{1}{p-1}}\,,
\end{equation}
if either
\par\noindent
(i) $1<q<\infty$, $q(p-1) \leq \sigma <\infty$ and
\begin{equation}\label{Ci}
\sup _{0< s < \m2 }\frac{s^{\frac{p-1}{\sigma}+\frac{1}{q'}}}{{\nu
_p} (s)} < \infty\,,
\end{equation}
or
\par\noindent
(ii) $1<q<\infty$, $0 < \sigma < q(p-1)$ and
\begin{equation}\label{Cii}
\int _0^{\m2}\bigg(\frac{s}{{\nu _p}(s)}\bigg)^{\frac{\sigma
q}{q(p-1)-\sigma}} ds < \infty\,,
\end{equation}
or
\par\noindent
(iii) $0<\sigma \le 1$, $q =\infty$ and
\begin{equation}\label{Ciii}
\int _0^{\m2}\Big(\frac{s}{{\nu _p}(s)}\Big)^{\frac{ \sigma}{p-1}}
ds < \infty\,.
\end{equation}
Moreover the constant $C$ in \eqref{Co} depends only on  $p$, $q$,
$\sigma$  and on the left-hand side either of \eqref{Ci}, or
\eqref{Cii}, or \eqref{Ciii}, respectively.
\end{theorem}
\par
\medskip

Theorem  \ref{teoA} is proved in  \cite{CM1}. Theorem \ref{teoC}
can be derived from Theorem  \ref{teoA}, via suitable weighted
Hardy type inequalities. In particular, a proof of cases $(i)$ and
$(ii)$ can be found in \cite[Theorem 4.1]{CM1}. Case $(iii)$
follows from case $(vi)$ of \cite[Theorem 4.1]{CM1}, via a
weighted Hardy type inequality for nonincreasing functions
\cite[Theorem 3.2 (b)]{HenigMaligranda}.
\par
\smallskip
We are now concerned with gradient estimates. A counterpart of
Theorem \ref{teoA} for $|\nabla u|$ is the content of the next
result.

\begin{theorem}\label{teo6}
Let $\o$, $p$ and $a$ be as in Theorem \ref{teo4bis}. Assume that
$f\in L^1(\o )\cap (V^{1,p}(\o ))'$ and fulfills \eqref{mv}. Let
$u$ be the weak solution to \eqref{1} satisfying $med(u)=0$. Then
\begin{equation}\label{confrgrad}
|\nabla u_\pm|^*(s)\le \left
(\frac{2}{s}\int_{\frac{s}{2}}^{\m2}\left (
\int_0^rf_{\pm}(\rho)\, d\rho   \right)^{p'}
d(-D\nu _p^{\frac{1}{1-p}})(r) \right )^\frac{1}{p}\qquad \hbox{for
$s\in (0, |\o|)$.}
\end{equation}
\end{theorem}

The proof of Theorem \ref{teo6} combines lower and upper estimates
for the integral of $|\nabla u|^{p-1}$ over the boundary of the
level sets of $u$. The relevant lower estimate involves the
isocapacitary function $\nu _p$. Given $u \in V^{1,p}(\o)$, we define
$\psi _u : [0, \infty ) \to [0, \infty )$ as
\begin{equation}\label{-519}
\psi _u (t) = \int _0^t \frac{d\tau}{\big(\int _{\{u=\tau \}}
|\nabla u|^{p-1}d \hh (x)\big)^{1/(p-1)}} \qquad \quad \hbox{for
$t\geq 0$}\,.
\end{equation}
As a consequence of \cite[Lemma 2.2.2/1]{Ma2}, one has that
\begin{equation}\label{-518}
C_p(\{u \geq t\},\{u>0\}) \leq  \psi _{u}(t)^{1-p} \qquad \quad
\hbox{for $t>0$}.
\end{equation}
Thus, if $u \in V^{1,p}(\o)$ and  fulfils \eqref{-520}, then on
estimating the infimum on the right-hand side of \eqref{-509} by
the choice
%
$E=\{u_\pm \geq t\}$ and $G=\{u_\pm
>0\}$, and on making use of \eqref{-518} applied with $u$ replaced by $u_+$ and
$u_-$, we deduce that
\begin{equation}\label{-522}
\nu _p(|\{u_\pm \geq t\}|) \leq \psi _{u_\pm }(t)^{1-p} \qquad \quad
\hbox{for $t > 0$}.
\end{equation}
\par
The upper estimate is contained in the following lemma from
\cite{CM1}, a version for Neumann problems of a result of
\cite{Ma3, Ta1, Ta3}.

\begin{lemma}\label{lemB}
Under the same assumptions as in Theorem \ref{teoA},
\begin{equation}\label{risB}
\int _{\{u_\pm = t\}}|\nabla u |^{p-1} d\hh (x) \leq \int _0^{\mu
_{u_\pm}(t)} f_\pm ^*(r)\,dr \qquad \quad \hbox{for a.e. $t>0$.}
\end{equation}
\end{lemma}

\bigskip

\noindent {\bf Proof of Theorem \ref{teo6}}. We shall prove
\eqref{confrgrad} for $u_+$, the proof for $u_-$ being analogous.
Consider the function $U:÷(0, \m2]\rightarrow [0, \infty)$ given
by
\begin{equation}\label{U}
U(s)=\int_{\{u_+\le u_+^*(s)\}} |\nabla u_+|^p\, dx, \qquad \qquad
\hbox{for $s\in (0, \m2]$.}
\end{equation}
Since $u\in W^{1,p}(\o)$, the function $u_+^*$ is locally
absolutely continuous (a.c., for short) in $(0, \m2)$ - see e.g.
\cite{CEG}, Lemma 6.6. The function
$$
(0, \infty) \ni t \mapsto \int_{\{u_+\le t\}}|\nabla u_+|^p\, dx
$$
is also locally a.c., inasmuch as, by coarea formula,
\begin{equation}\label{coarea}
\int_{\{u_+\le t\}}|\nabla u_+|^p\, dx=\int_0^t\int_{\{u_+=
\tau\}}   |\nabla u_+|^{p-1}d\hh (x) d\tau, \qquad \hbox{for }
t>0.
\end{equation}

Thus, $U$ is locally a.c., for it is the composition of monotone
a.c. functions, and by \eqref{coarea}
\begin{equation}\label{coareabis}
U'(s)=-{u_+^*}'(s)\int_{\{u_+= u_+^*(s)\}} |\nabla u_+|^{p-1}d\hh
(x), \quad \hbox{for a.e. $s\in(0, \m2)$}.
\end{equation}
Similarly, the function
$$
(0, \m2) \ni s \mapsto \psi_{u_+}(u_+^*(s)),
$$
where $\psi_{u_+}$  is defined as in
 \eqref{-519},
is locally a.c., and
\begin{equation}\label{derpsi}
\frac{d}{ds}\left ( \psi_{u_+}(u_+^*(s)) \right
)=\frac{{u_+^*}'(s)}{\int_{\{u_+= u_+^*(s)\}} |\nabla
u_+|^{p-1}d\hh (x)} \qquad \hbox{for a.e. $s\in(0, \m2)$}.
\end{equation}
Let us set
\begin{equation*}
 W(s) = \frac{d}{ds}\left ( \psi_{u_+}(u_+^*(s))\right) \quad \qquad
\hbox{for a.e. $s\in(0, \m2)$}.
\end{equation*}
 From
\eqref{coareabis}, \eqref{derpsi} and \eqref{risB}, we obtain that
\begin{equation}\label{611}
-U'(s)\le W(s) \left (\int_0^sf^*_+(r)\, dr  \right )^{p'}, \qquad
\hbox{for a.e. $s\in(0, \m2)$}.
\end{equation}
Note that in deriving \eqref{611}  we have made use of the fact
that $\mu _{u_+}(u_+^*(s))=s$ if $s$ does not belong to any
interval where $u_+^*$ is constant, and that ${u_+^*}'=0$ in any
such interval. Since $\psi_{u_+}(u_+^*(\m2))=\psi_{u_+}(0)=0$,
from \eqref{-522} we obtain that
\begin{equation}\label{613}
\int_s^{\m2}W(r)\, dr=\psi_{u_+}(u_+^*(s))\le
\nu _p^{\frac{1}{1-p}}(s)=\int_{s}^{\m2}d(-D\nu _p^{\frac{1}{1-p}})(r),
\qquad \hbox{for } s\in(0, \m2)\,.
\end{equation}
Owing to Hardy's lemma (see e.g. \cite[Chapter 2, Proposition
3.6]{BS}), inequality \eqref{613} entails that
\begin{equation}\label{614}
\int_0^{\m2}\phi(r) W(r) dr \le
\int_0^{\m2}\phi(r)d(-D\nu _p^{\frac{1}{1-p}})(r)
\end{equation}
for every non-decreasing function $\phi:÷(0, \m2)\rightarrow [0,
\infty)$.
 In particular, fixed any such function $\phi$, we have
 that
\begin{equation}\label{615}
\int_0^{\m2}\phi(r)  \left (\int_0^rf^*_+(\rho)\, d\rho  \right
)^{p'}W(r) dr\le \int_0^{\m2}\phi(r)  \left (\int_0^rf^*_+(\rho)\,
d\rho  \right )^{p'} d(-D\nu _p^{\frac{1}{1-p}})(r).
\end{equation}
 Coupling \eqref{611} and
\eqref{615} yields
\begin{equation}\label{616}
\int_0^{\m2}-U'(r)\phi(r) dr\le \int_0^{\m2}\phi(r)  \left
(\int_0^rf^*_+(\rho)\, d\rho  \right )^{p'}
d(-D\nu _p^{\frac{1}{1-p}})(r).
\end{equation}
Next note that
\begin{equation}\label{617}
\int_{s}^{\m2}-U'(r)dr=U(s)=\int_{\{u_+\le u^*(s)\}}|\nabla u_+|^p dx\ge \int_{s}^{\m2}|\nabla u_+|^*(r)^pdr
\end{equation}
for $s\in(0, \m2)$, where the inequality follows from the first
inequality in \eqref{HL} and from the inequality $|\{0\le u_+\le
u_+^*(s) \}|\ge \m2-s.$ Inequality \eqref{617}, via Hardy's lemma
again, ensures that
\begin{equation}\label{621}
\int_0^{\m2}|\nabla u_+|^*(r)^p\phi(r)dr\le
\int_0^{\m2}-U'(r)\phi(r)dr\,.
\end{equation}
Fixed any $s\in (0,\m2)$, we infer from \eqref{616} and
\eqref{621} that
\begin{equation}\label{622}
|\nabla u_\pm|^*(s)^p\int_0^s\phi(r)dr\le \int_0^{\m2}
\phi(r)\left ( \int_0^rf_{\pm}(\rho)\, d\rho   \right)^{p'}
d(-D\nu _p^{\frac{1}{1-p}})(r).
\end{equation}
 Inequality \eqref{confrgrad} follows
from \eqref{622} on choosing $\phi =\chi_{[s/2,\m2]}$. \qed
\medskip

Estimates for Lebesgue norms of $|\nabla u|$ are provided by the
next result.

\begin{theorem}\label{teo5ter}
Let $\o$, $p$ and $a$ be as in Theorem \ref{teo4bis}. Assume that
 $f\in L^q(\o)\cap (V^{1,p}(\o ))'$ for some $q\in [1, \infty]$ and fulfills \eqref{mv}.
Let $u$ be a weak solution to problem \eqref{1}. Let $0<\sigma\le
p$. Then there exists a constant $C$ such that
\begin{equation}\label{stima5ter}
\|\nabla u \|_{L^\sigma (\o )} \leq C \|f\|_{L^q(\o
)}^{\frac{1}{p-1}}\,,
\end{equation}
if either
\par\noindent
(i) $q>1$, $q(p-1) \leq \sigma $ and
\begin{equation}\label{5teri}
\sup _{0< s < \M2
}\frac{s^{1+\frac{p(p-1)}{\sigma}-\frac{p}{q}}}{{\nu _p}(s)} <
\infty\,,
\end{equation}
or
\par\noindent
(ii) $1<q<\infty$, $0 < \sigma < q(p-1)$ and
\begin{equation}\label{5terii}
\int _0^{\m2}\bigg(\frac{s}{{\nu _p}(s)}\bigg)^{\frac{\sigma
q}{p[q(p-1)-\sigma]}} ds < \infty\,,
\end{equation}
or
\par\noindent
(iii)  $q =\infty$ and
\begin{equation}\label{5teriii}
\int _0^{\m2}\Big(\frac{s}{{\nu _p}(s)}\Big)^{\frac{
\sigma}{p(p-1)}} ds < \infty\,,
\end{equation}
or
\par\noindent
(iv)  $q =1$ and
\begin{equation}\label{5teriv}
\int _0^{\m2}\Big(\frac{s}{{\nu _p}(s)}\Big)^{\frac{
\sigma}{p(p-1)}} \frac{ds}{s^\frac{\sigma}{p-1}} < \infty\,.
\end{equation}
Moreover the constant $C$ in \eqref{stima5ter} depends only on
$p$, $q$, $\sigma$  and on the left-hand side either of
\eqref{5teri}, or \eqref{5terii}, or \eqref{5teriii} or
\eqref{5teriv}, respectively.
\end{theorem}
\medskip

Cases \emph{(i)--(iii)} of Theorem \ref{teo5ter} are proved in
\cite[Theorem 5.1]{CM1}; an alternative proof can be given by an
argument analogous to that of Theorem \ref{teo5}, Section
\ref{sec4}. Case  \emph{(iv)} is a straightforward consequence of
the following proposition.

\begin{proposition}\label{prop7}
Let $\o$, $p$ and $a$ be as in Theorem \ref{teo4bis}. Assume that
$f\in L^1(\o )\cap (V^{1,p}(\o ))'$ and fulfills \eqref{mv}. Let
$u$ be a weak solution to \eqref{1}.
%
Let $\omega _p: (0, |\o|) \rightarrow [0,\infty)$ be the
function defined by
\begin{equation}\label{540'}
\omega _p(s)=(s{\nu _p}^{\frac{1}{p-1}}(s/2))^\frac{1}{p},
\qquad \hbox{for }s\in(0, |\o |).
\end{equation}
Then there exists a constant $C=C(p, n)$ such that
\begin{equation}\label{540}
\|\nabla u \|_{M_{\omega _p}(\o )} \leq C \|f\|_{L^1(\o
)}^{\frac{1}{p-1}}\,,
\end{equation}
where $M_{\omega _p}(\o )$ is the Marcinkiewicz space
defined as in \eqref{normmarc}
\end{proposition}

\noindent{\bf Proof}. If $u$ is normalized in such a way that
${\rm med } (u)=0$, by estimate \eqref{confrgrad} one gets that
\begin{equation*}
|\nabla u_\pm|^*(s)\le  \|
f_{\pm}\|_{L^1(\o)}^{\frac{1}{p-1}}\left
(\frac{2}{s}\int_{\frac{s}{2}}^{\m2}  d(-D\nu _p^{\frac{1}{1-p}})(r)
\right )^\frac{1}{p}\le 2^\frac{1}{p}\|
f_{\pm}\|_{L^1(\o)}^{\frac{1}{p-1}} \left (  \frac{1}{s}
\nu _p^{\frac{1}{1-p}}(s/2)  \right)^\frac{1}{p}
\end{equation*}
for $s\in (0, |\o|)$. Inequality \eqref{540} follows.
\qed

Let us note that Theorem \ref{teo6} can also be used to provide a
further alternate proof of Cases (i)--(iii) of Theorem
\ref{teo5ter} when $\sigma < p$. In fact, these cases are special
instances of Theorem \ref{teo9} below, dealing with a priori
estimates for Lorentz norms of the gradient. Theorem \ref{teo9} in
turn rests upon the following corollary of Theorem \ref{teo6}.

\begin{corollary}\label{cor8}
Let $\o$, $p$ and $a$ be as in Theorem \ref{teo4bis}. Let $X(\o)$
be an r.i. space and let $f\in X(\o)\cap (V^{1,p}(\o ))'$. Let $u$
be a weak solution to problem \eqref{1}. Assume that $Y(\o)$ is an
r.i. space such that
\begin{equation}\label{801}
\left \|     \left (\frac{1}{s}\int_s^{|\o|}\left (
\int_0^r\phi(\rho)\, d\rho   \right)^{p'}
d(-D\nu _p^{\frac{1}{1-p}})(r) \right )^\frac{1}{p} \right
\|_{\overline Y(0, |\o|)}\le C \|\phi\|^\frac{1}{p-1}_{\overline
X(0, |\o|)},
\end{equation}
for some constant $C$ and every nonnegative and non-increasing
function $\phi\in \overline X(0, |\o|)$. Then there exists a
constant $C_1=C_1(C)$ such that
\begin{equation}\label{802}
\|\nabla u\|_{Y(\o)}\le C_1\|f\|^\frac{1}{p-1}_{ X(\o )}.
\end{equation}
\end{corollary}
\noindent{\bf Proof.} Inequality \eqref{802} immediatly follows
from \eqref{confrgrad}, \eqref{801} and the fact that the dilation
operator  $H$ defined on any  function $\phi \in \Mc (0, |\o|)$ by
$$
H\phi (s)=\phi (s/2), \qquad \hbox{for }s\in (0, |\o|),
$$
is bounded in any r.i. space on $(0, |\o|)$ (see e.g.
\cite[Chapter 3, Proposition 5.11]{BS}). \qed
\begin{remark}\label{maximal}
{\rm If $X(\o)$ is such that the Hardy type inequality
\begin{equation}\label{802bis}
\left \| \frac{1}{s}\int_0^s\phi(r)dr \right \|_{\overline X(0,
|\o|)}\le C_2 \|
 \phi
 \|_{\overline X(0, |\o|)}
\end{equation}
holds every nonnegative and non-increasing function $\phi\in
\overline X(0, |\o|)$ and for some constant $C_2$, and
\begin{equation}\label{803}
\left \|     \left (\frac{1}{s}\int_s^{|\o|} \phi(r)^{p'} r^{p'}
d(-D\nu _p^{\frac{1}{1-p}})(r) \right )^\frac{1}{p} \right
\|_{\overline Y(0, |\o|)}\le C_3 \|\phi\|^\frac{1}{p-1}_{\overline
X(0, |\o|)},
\end{equation}
for some constant $C_3$ and every $\phi$ as above, then
\eqref{802} holds with $C_1=C_1(C_2,C_3)$. \noindent Indeed, if
\eqref{803} is in force, then \eqref{801}  holds with $\phi$
replaced by $\frac{1}{s}\int_0^s \phi(r)dr$ on the right-hand
side. Inequality \eqref{802} then follows via \eqref{802bis}. }
\end{remark}
%

\begin{theorem}\label{teo9}
Let $\o$, $p$ and $a$  be as in Theorem \ref{teo4bis}.  Let
$0<\sigma<p$, $1<q< \infty $, $0<\gamma,\varrho<\infty$. Let $f\in
L^{q, \frac{\gamma}{p-1}}(\o)\cap (V^{1,p}(\o ))'$ and let $u$ be
a weak solution to problem \eqref{1}.  Then there exists a
constant $C$ such that
\begin{equation}\label{623}
\|\nabla u \|_{L^{\sigma,\varrho} (\o )} \leq C \|f\|_{L^{q,
\frac{\gamma}{p-1}}(\o )}^{\frac{1}{p-1}}\,
\end{equation}
if either
\par\noindent
(i) $\gamma \le \varrho$ and
\begin{equation}\label{624}
\sup _{0< s < \M2
}\frac{s^{1+\frac{p(p-1)}{\sigma}-\frac{p}{q}}}{{\nu _p}(s)} <
\infty\,,
\end{equation}
or
\par\noindent
(ii) $\gamma>\varrho$ and
\begin{equation}\label{625}
\int
_0^{\m2}\bigg(\frac{s^{1+\frac{p(p-1)}{\sigma}-\frac{p}{q}}}{{\nu
_p} (s)}\bigg)^{\frac{\varrho\gamma}{p(\gamma -\varrho)(p-1)}}
\frac{ds}{s} < \infty\,.
\end{equation}
Moreover the constant $C$ in \eqref{623} depends only on $p$, $q$,
$\sigma$, $\varrho$, $\gamma$ and on the left-hand side either of
\eqref{624} or \eqref{625}, respectively.
\end{theorem}
\medskip

The proof of Theorem \ref{teo9} relies  upon Corollary \ref{cor8}
and on a characterization of weighted one-dimensional Hardy-type
inequalities for non-increasing functions established in
\cite{Gold}. The arguments to be used are similar to those
exploited in the proof of \cite[Theorem 4.1]{CM1}. The details are
omitted for brevity.

\subsection{Proof of Theorem \ref{teo4bis}}\label{subsec3.2}

A key step in our proof of Theorem \ref{teo4bis} is the following
uniform integrability result for the gradient of weak solutions to
\eqref{1} with  $f\in (V^{1,p}(\o ))'$, which relies upon Theorem
\ref{teo6}.

\begin{lemma}\label{lemnuovo}
Let $\o$, $p$ and $a$ be as in Theorem \ref{teo4bis}. Assume that
 $f\in L^q(\o)\cap (V^{1,p}(\o ))'$ for some $q\in [1, \infty]$ and fulfills \eqref{mv}.
  Then there exists a
function $\varsigma : (0, \infty ) \to [0, \infty )$, depending on
$\Omega$, $p$ and $q$, satisfying
\begin{equation}\label{n-2}
\lim _{s \to 0^+} \varsigma (s) = 0\,,
\end{equation}
and such that, if $u$ is a weak solution to \eqref{1} satisfying
\eqref{-520}, then
\begin{equation}\label{n-1}
\int _F |\nabla u_{\pm}|^{p-1} dx \leq \varsigma (|F|)\,
\|f_{\pm}\|_{L^{q}(\Omega )}
\end{equation}
for every measurable set $F \subset \Omega$.
\end{lemma}
{\bf Proof}. By the Hardy-Littlewood inequality \eqref{HL} and
Theorem \ref{teo6}, we have that
\begin{align}\label{n1}
\int _F |\nabla u_{\pm}|^{p-1} dx & \leq \int _0^{|F|} |\nabla
u_{\pm}|^*(s)^{p-1} ds
\\
\nonumber & \leq 2^{1/p}  \int _0^{|F|/2}\bigg(\frac 1s \int
_s^{\m2} \bigg(\int _0^r f_{\pm}^*(\rho )d\rho \bigg)^{p'}
 d(-D\nu _p^{\frac{1}{1-p}})(r) \bigg )^{\frac{1}{p'}}ds
\\
\nonumber &\leq
 2^{1/p}  \int _0^{|F|/2}\bigg(\frac 1s \int
_s^{|F|/2} \bigg(\int _0^r f_{\pm}^*(\rho )d\rho \bigg)^{p'}
 d(-D\nu _p^{\frac{1}{1-p}})(r) \bigg)^\frac{1}{p'}ds
\\ \nonumber
\quad &+ p\, 2^{1/p} \Big(\frac{|F|}{2}\Big)^{1/p} \bigg(\int
_{|F|/2}^{\m2} \bigg(\int _0^r f_{\pm}^*(\rho )d\rho \bigg)^{p'}
 d(-D\nu _p^{\frac{1}{1-p}})(r) \bigg )^\frac{1}{p'}\,.
\end{align}
Assume first that $1<q\leq \infty$ and \eqref{401} is in force.
Let us preliminarily observe that
\begin{equation}\label{n2}
\lim _{s \to 0^+} \Big(\frac s{{\nu _p}(s)}\Big)^{\frac {q'}{p}}s
=0\,,
\end{equation}
since
\begin{equation*}
\int _0^s \Big(\frac r{{\nu _p}(r)}\Big)^{\frac {q'}{p}}dr \geq
\frac 1{{\nu _p}(s)^{\frac {q'}{p}}}\int _0^s r^{\frac {q'}{p}} dr
= \frac p{q'+p}\frac{s^{\frac {q'}{p}+1}}{{\nu _p}(s)^{\frac
{q'}{p}}} \quad \hbox{for $s \in (0, \m2 )$\,.}
\end{equation*}
Consider the second addend on the rightmost side of \eqref{n1}. We
claim that there exists a function $\kappa : (0, \infty ) \to [0,
\infty )$ such that
\begin{equation}\label{n4}
\lim _{s\to 0^+} s^{1/p} \kappa (s) =0\,,
\end{equation}
and
\begin{equation}\label{n5}
\bigg(\int _{s}^{\m2} \bigg(\int _0^r f_{\pm}(\rho )d\rho
\bigg)^{p'}
 d(-D\nu _p^{\frac{1}{1-p}})(r) \bigg )^\frac{1}{p'} \leq \kappa (s)
\|f_{\pm}^*\|_{L^q(0,\m2 )}\,.
\end{equation}
To verify this claim, assume first that $p' \geq q$. By a weighted
Hardy inequality \cite[Section 1.3]{Ma2}, inequality \eqref{n5}
holds with
\begin{equation}\label{n6}
\kappa (s) = 
C \sup _{s\leq r \leq \m2} {\nu _p}(r) ^{-1/p} r^{\frac 1{q'}}\,,
\end{equation}
for some constant $C=C(p,q)$. Moreover, $\kappa$ fulfils
\eqref{n4}, since
\begin{equation}\label{n7}
\lim _{s \to 0^+} s^{1/p}\sup _{s\leq r \leq \m2} {\nu _p}(r)
^{-1/p} r^{\frac 1{q'}}=0\,.
\end{equation}
Indeed, equation \eqref{n7} holds trivially if $\displaystyle{\sup
_{0< r \leq \m2} }{\nu _p}(r) ^{-1/p} r^{\frac 1{q'}}< \infty$. If
this is not the case, then for each $s \in (0, \m2)$ define
$$r(s) = \inf \Big\{r \in [s, \m2]: 2{\nu _p}(r) ^{-1/p}
r^{\frac 1{q'}}\geq \sup _{s\leq \rho \leq \m2} {\nu _p}(\rho)
^{-1/p} \rho ^{\frac 1{q'}}\Big\}$$ and observe that the function
$r(s)$ converges monotonically to $0$ as $s$ goes to $0$, and that
\begin{align}\label{n8}
\lim _{s \to 0^+} \Big(s^{1/p}\sup _{s\leq r \leq \m2} {\nu _p}(r)
^{-1/p} r^{\frac 1{q'}}\Big) & \leq 2 \lim _{s \to 0^+} s^{1/p}
{\nu _p} (r(s)) ^{-1/p} r(s)^{\frac 1{q'}} \\ \nonumber & \leq
2\lim _{s \to 0^+} \Big(\frac {r(s)}{{\nu _p}(r(s))}\Big)^{\frac
{1}{p}}r(s)^{1/q'}=0\,,
\end{align}
by \eqref{n2}.
\par\noindent
Consider next the case when $p' <q$. An appropriate weighted Hardy
inequality \cite[Section 1.3]{Ma2} now tells us that inequality
\eqref{n5} holds with
\begin{equation}\label{n9}
\kappa(s) = C \bigg(\int _0^{\m2} \bigg(r^{\frac 1{p-1}}\int
_r^{\m2} \chi _{(s, \m2 )}(\rho
)d(-D\nu _p^{\frac{1}{1-p}})(\rho)\bigg)^{\frac q{q-p'}}
dr\bigg)^{\frac {q-p'}{qp'}}\,,
\end{equation}
for some  constant $C=C(p,q)$. Here, the exponents $\frac q{q-p'}$
and $\frac {q-p'}{qp'}$ are replaced by $1$ and $\frac 1{p'}$,
respectively, when $q=\infty$. We have that
\begin{align}\label{n10}
\kappa (s) &= C\bigg(\int _0^s \big(r^{\frac 1{p-1}}\nu
_p(s)^{\frac 1{1-p}}\big)^{\frac q{q-p'}} dr + \int
_s^{\m2}\big(r^{\frac 1{p-1}}\nu _p(r)^{\frac
1{1-p}}\big)^{\frac q{q-p'}} dr\bigg)^{\frac {q-p'}{qp'}} \\
\nonumber & \leq C \frac {s^{1/q'}}{{\nu _p}(s)^{1/p}}+ C
\bigg(\int _s^{\m2}\bigg(\frac{r}{\nu _p(r)}\bigg)^{\frac
{q}{(p-1)(q-p')}}dr\bigg)^{\frac {q-p'}{qp'}}.
\end{align}
 Thus,
\begin{align}\label{n11}
s^{1/p} \kappa (s) & \leq C \Big(\frac s{{\nu _p}(s)}\Big)^{\frac
{1}{p}}s^{1/q'} + C s^{1/p}\bigg(\int _s^{\m2}\bigg(\frac{r}{\nu
_p(r)}\bigg)^{\frac {q}{(p-1)(q-p')}}dr\bigg)^{\frac {q-p'}{qp'}}
\\ \nonumber
 &\leq C
\Big(\frac s{{\nu _p}(s)}\Big)^{\frac {1}{p}}s^{1/q'} + C s^{1/p}
\sup _{s\leq r \leq \m2} \Big(\frac{r}{{\nu _p}(r)}\Big)^{\frac
{q'}{p^2}}\bigg(\int _s^{\m2}\bigg(\frac{r}{\nu
_p(r)}\bigg)^{\frac {q'}{p}}dr\bigg)^{\frac {q-p'}{qp'}}\,,
\end{align}
and hence $\kappa (s)$ fulfils \eqref{n4} also in this case, by
\eqref{n2}, by \eqref{401} and by the fact that
$$\lim _{s \to
0^+}s^{1/p} \Big(\sup _{s\leq r \leq \m2} \Big(\frac{r}{{\nu _p}
(r)}\Big)^{\frac {q'}{p^2}}\Big) =0,$$
as an analogous argument as
in the proof of    \eqref{n8} shows.
\par\noindent
We have thus proved that
\begin{equation}\label{n12}
s^{1/p}\bigg(\int _{s}^{\m2} \bigg(\int _0^r f_{\pm}(\rho )d\rho
\bigg)^{p'}
 d(-D\nu _p^{\frac{1}{1-p}})(r) \bigg )^\frac{1}{p'}\leq \varsigma (|E|)\,
\|f_{\pm}^*\|_{L^{q}(0, \m2 )} \qquad \hbox{for $s \in (0, \m2
)$},
\end{equation}
for some  function $\varsigma$ as in the statement.
\par
Let us now take into account the first addend on the rightmost
side of \eqref{n1}. We shall show that
\begin{multline}\label{n13}
\int _0^{s}r^{-1/p'}\bigg( \int _r^s  f_{\pm}^{**}(\rho )^{p'}\rho
^{p'}
 d(-D\nu _p^{\frac{1}{1-p}})(\rho ) \bigg)^\frac{1}{p'}dr \\ \leq C
\bigg(\int _0^{s}\bigg(\frac{r}{\nu _p(r)}\bigg)^{\frac {q'}{p}} dr
\bigg)^\frac{1}{q'} \|f_{\pm}^{**}\|_{L^{q}(0,s)} \qquad \hbox{for
$s \in (0, \m2 )$},
\end{multline}
for some constant $C=C(p,q)$. It suffices to establish \eqref{n13}
for some fixed number $s$, say $1$, since the general case then
follows by scaling. As a consequence of \cite[Theorem 1.1 and
Remark 1.4]{Gold}, the inequality
\begin{equation}\label{n14}
\int _0^{1}r^{-1/p'}\bigg( \int _r^1  \phi(\rho )^{p'}\rho ^{p'}
 d(-D\nu _p^{\frac{1}{1-p}})(\rho ) \bigg)^\frac{1}{p'}dr \leq C
 \|\phi\|_{L^{q}(0,1)}
\end{equation}
holds for every nonnegative non-increasing function $\phi$ in
$(0,1)$ if
\begin{equation}\label{n15}
\bigg(\int _0^1 \bigg(\int _0^r \bigg(\rho ^{p'} {\nu _p} (\rho
)^{\frac 1{1-p}} + \int _\rho ^r \bigg( \frac {\theta}{{\nu _p}
(\theta )}\bigg)^{\frac 1{p-1}} d\theta \bigg)^{\frac 1{p'}} \rho
^{-\frac 1{p'}} d\rho\bigg)^{q'}r^{-q'}dr\bigg)^{\frac 1{q'}} <
\infty\,.
\end{equation}
Moreover, the constant $C$ on the right-hand side of \eqref{n14}
does not exceed the integral on the left-hand side of \eqref{n15}
(up to a multiplicative constant depending on $p$ and $q$). Thus,
inequality \eqref{n13} will follow if we show that
\begin{equation}\label{n16}
\int _0^1 \bigg(\int _0^r \bigg(\rho ^{p'} {\nu _p} (\rho )^{\frac
1{1-p}} + \int _\rho^r \bigg( \frac {\theta}{{\nu _p} (\theta
)}\bigg)^{\frac 1{p-1}} d\theta \bigg)^{\frac 1{p'}} \rho ^{-\frac
1{p'}} d\rho\bigg)^{q'}r^{-q'}dr \leq C \int
_0^{1}\bigg(\frac{r}{\nu _p(r)}\bigg)^{\frac {q'}{p}} dr\,,
\end{equation}
for some constant $C=C(p,q)$. The standard Hardy inequality
entails that
\begin{equation}\label{n17}
\int _0^1 \bigg(\int _0^r \bigg(\frac {\rho}{\nu _p(\rho
)}\bigg)^{\frac 1p} d\rho\bigg) ^{q'} r^{-q'}dr \leq C
 \int
_0^{1}\bigg(\frac{r}{\nu _p(r)}\bigg)^{\frac {q'}{p}} dr\,,
\end{equation}
for some constant $C=C(q)$. Thus, it only remains to prove that
\begin{equation}\label{n18}
\int _0^1 \bigg(\int _0^r \bigg(\int _\rho^r \bigg( \frac
{\theta}{{\nu _p} (\theta )}\bigg)^{\frac 1{p-1}} d\theta
\bigg)^{\frac 1{p'}} \rho ^{-\frac 1{p'}}
d\rho\bigg)^{q'}r^{-q'}dr \leq C \int _0^{1}\bigg(\frac{r}{\nu
_p(r)}\bigg)^{\frac {q'}{p}} dr\,,
\end{equation}
for some constant $C=C(p,q)$. Consider first the case when $p<q$.
By a Hardy type inequality again, we have that
\begin{align}\label{n19}
\int _0^1  \bigg(\int _0^r \bigg(&\int _\rho^r \bigg( \frac
{\theta}{{\nu _p} (\theta )}\bigg)^{\frac 1{p-1}} d\theta
\bigg)^{\frac 1{p'}} \rho ^{-\frac 1{p'}}
d\rho\bigg)^{q'}r^{-q'}dr \\ \nonumber & \leq \int _0^1 \bigg(\int
_0^r \bigg(\int _\rho ^1 \bigg( \frac {\theta}{{\nu _p} (\theta
)}\bigg)^{\frac 1{p-1}} d\theta \bigg)^{\frac 1{p'}} \rho ^{-\frac
1{p'}} d\rho\bigg)^{q'}r^{-q'}dr
\\ \nonumber & \leq C \int
_0^1 \bigg(\int _r ^1 \bigg( \frac {\theta}{{\nu _p} (\theta
)}\bigg)^{\frac 1{p-1}} d\theta \bigg)^{\frac {q'}{p'}}
r^{-\frac{q'}{p'}}dr\,.
\end{align}
On the other hand, since ${\nu _p}$ is a non-increasing function,
by \cite[Theorem 1.1 and Remark 1.4]{Gold}, the right-hand side of
\eqref{n19} does not exceed the right-hand side of \eqref{n18},
and hence \eqref{n18} follows. Assume now that $p \geq q$. Then
\begin{align}\label{n20}
\int _0^1  \bigg(\int _0^r \bigg(&\int _\rho^r \bigg( \frac
{\theta}{{\nu _p} (\theta )}\bigg)^{\frac 1{p-1}} d\theta
\bigg)^{\frac 1{p'}} \rho ^{-\frac 1{p'}}
d\rho\bigg)^{q'}r^{-q'}dr  \\ \nonumber &\leq
 \int _0^1
\bigg(\int _0^r \bigg( \frac {\theta}{{\nu _p} (\theta
)}\bigg)^{\frac 1{p-1}} d\theta \bigg)^{\frac {q'}{p'}}\bigg(\int
_0^r \rho ^{-\frac 1{p'}} d\rho\bigg)^{q'}r^{-q'}dr
 \\ \nonumber &\leq C
\int _0^{1}\bigg(\frac{r}{\nu _p(r)}\bigg)^{\frac {q'}{p}} dr\,,
\end{align}
for some constant $C=C(p,q)$, where the last inequality holds by
the Hardy inequality. Inequality \eqref{n18} is established also
in this case. Thus, inequality \eqref{n16}, and hence \eqref{n13},
is fully proved. Combining \eqref{n1}, \eqref{n12} and
\eqref{n13}, and making use of the fact that
$$\|f_{\pm}^{**}\|_{L^{q}(0,s)} \leq C
\|f_{\pm}^{*}\|_{L^{q}(0,s)}$$ for some constant $C=C(q)$, by the
Hardy inequality, conclude the proof in the case when $1<q\leq
\infty$.
\par
Let us finally focus on the case when $q=1$. The counterpart of
\eqref{n2} is now
\begin{equation}\label{n21}
\lim _{s\to 0^+} \frac{s}{{\nu _p} (s) } =0\,.
\end{equation}
Inequality \eqref{n5} holds with
\begin{equation}\label{n22}
\kappa (s) = C {\nu _p} (s)^{-1/p}\,,
\end{equation}
and hence, by \eqref{n21}, the function $\kappa$ fulfils
\eqref{n4}. Inequality \eqref{n12} is thus established. On the
other hand,
\begin{align}\label{n23}
\int _0^{s}r^{-1/p'}& \bigg( \int _r^s  f_{\pm}^{**}(\rho
)^{p'}\rho ^{p'}
 d(-D\nu _p^{\frac{1}{1-p}})(\rho ) \bigg)^\frac{1}{p'}dr \\ \nonumber
&\leq \|f_{\pm}^*\|_{L^1(0, \m2 )} \int _0^{s}r^{-1/p'} \bigg(
\int _r^s
 d(-D\nu _p^{\frac{1}{1-p}})(\rho ) \bigg)^\frac{1}{p'}dr
\\ \nonumber
&\leq \|f_{\pm}^*\|_{L^1(0, \m2 )} \int _0^{s}\bigg(\frac r{\nu
_p(r)}\bigg)^{1/p} \frac{ dr}{r}\,.
\end{align}
The conclusion follows via \eqref{n1}, \eqref{n12} and
\eqref{n23}. \qed
\smallskip
\par
We are now ready to prove Theorem \ref{teo4bis}.
\smallskip
\par\noindent
{\bf Proof of Theorem \ref{teo4bis}} Our assumptions ensure that a
sequence $\{ f_k \}\subset L^q(\o)\cap (V^{1,p}(\o ))'$ exists
such that
\begin{equation}\label{404}
f_k \rightarrow f \qquad \hbox{in }L^q(\o)
\end{equation}
and
\begin{equation}\label{405}
\int_\o f_k \,dx  =0 \qquad \hbox{for } k\in \N.
\end{equation}
Indeed, if $1\le q<\infty$, any sequence $\{ f_k \}$ of continuous
compactly supported functions fulfilling \eqref{404} and
\eqref{405} does the job; when $q=\infty$, it suffices to take
$f_k=f$ for $k\in \N$, since $V^{1,p}(\o) \to L^1(\o)$ provided
that \eqref{402} is in force, by \eqref{2002'}. We may also
clearly assume that
\begin{equation}\label{406}
\|f_k\|_{L^q(\o)}\le 2\|f\|_{L^q(\o)}, \qquad \hbox{for } k\in \N.
\end{equation}
By Proposition \ref{weak}, for each $k\in \N$ there exists an
unique weak solution $u_k\in V^{1,p}(\o)$ to the problem
\begin{equation}\label{407}
\begin{cases}
- {\rm div} (a(x,\nabla u_k)) = f_k(x)  & {\rm in}\,\,\, \o \\
 a(x,\nabla u_k)\cdot {\bf  n} =0 &
{\rm on}\,\,\,
\partial \o \,
\end{cases}
\end{equation}
fulfilling
\begin{equation}\label{408}
{\rm med } (u_k)=0.
\end{equation}
Hence,
\begin{equation}\label{408bis}
\int_\o a(x, \nabla u_k)\cdot \nabla \Phi \, dx=\int_\o f_k\Phi\,
dx\, ,
\end{equation}
for every $\Phi \in V^{1, p}(\o)$.

We split the proof of the existence of an approximable solution to
\eqref{1} in steps. The outline of the argument is related to that
of \cite{BBGGPV, DMOP}. \par\noindent

\medskip
\noindent{\bf Step 1.} \emph{There exists a measurable function
$u: \Omega \to \R$ such that}
\begin{equation}\label{408ter}
u_k\rightarrow u \qquad \hbox{a.e. in }\o,
\end{equation}
\emph{up to subsequences. Hence, property (ii) of the definition
of approximable solution holds}.

Given any $t, \tau>0$, one has that
\begin{equation}\label{409}
|\{|u_k-u_m|>\tau  \}|\le |\{| u_k|>t \}|+|\{| u_m|>t \}|+
|\{|\T(u_k)-\T(u_m)|>\tau  \}|,
\end{equation}
for $k,m\in \N$. By \eqref{A} and \eqref{406},
\begin{equation}\label{410}
(u_k)_\pm ^* (s) \leq
\nu _p(s)^{\frac{1}{1-p}}\|(f_k)_\pm\|_{L^1(\o)}^{\frac{1}{p-1}}\le
2^{\frac{1}{p-1}}|\o|^{\frac{1}{q'(p-1)}}\nu _p(s)^{\frac{1}{1-p}}\|f\|_{L^q(\o)}^{\frac{1}{p-1}}
\quad \hbox{for }s\in(0, |\o |/2),
\end{equation}
and for $k\in \N$, whence
\begin{equation}\label{411}
\mu_{(u_k)_\pm}(t)\le
\nu _p^{-1}\left(\frac{2|\o|^{\frac{1}{q'}}\|f\|_{L^q(\o)}}{t^{p-1}}\right),
\qquad \hbox{for }t>0,
\end{equation}
and for $k\in \N$. Here, $\nu _p^{-1}$ denotes the generalized
left-continuous inverse of $\nu _p$. Thus, fixed any $\varepsilon>0$,
the number $t$ can be chosen so large that
\begin{equation}\label{412}
 |\{| u_k|>t \}|<\varepsilon  \quad \hbox{and} \quad  |\{| u_m|>t \}|<\varepsilon.
\end{equation}
Next, fix any smooth open set $\o_\varepsilon\subset \subset \o$ such that
\begin{equation}\label{413}
|\o \setminus \o_\varepsilon|<\varepsilon.
\end{equation}
On choosing $\Phi=\T(u_k)$ in \eqref{408bis} and making use of
\eqref{406} we obtain that
\begin{equation}\label{414}
\int_\o|\nabla \T(u_k)|^p\, dx =\int_{\{|u_k|<t\}}|\nabla u_k|^p\, dx \le \int_{\{|u_k|<t\}}a(x, \nabla u_k)\cdot
\nabla u_k\, dx\le 2t|\o|^{\frac{1}{q'}}\|f\|_{L^q(\o)},
\end{equation}
for $k\in \N$. In particular the sequence $\{\T(u_k)\}$  is
bounded in $W^{1,p}(\o_\varepsilon)$. By the compact embedding of
$W^{1,p}(\o_\varepsilon)$ into $L^p(\o_\varepsilon)$, $\T(u_k)$
converges  (up to subsequences) to some function in
$L^p(\o_\varepsilon)$. In particular, $\{\T(u_k)\}$ is a Cauchy
sequence in measure in $\o_\varepsilon$. Thus,
\begin{equation}\label{416}
|\{ |\T(u_k)-\T(u_m)|>\tau  \}|\le |\o\setminus
\o_\varepsilon|+|\o _\varepsilon \cap \{ |\T(u_k)-\T(u_m)|>\tau
\}|<2\varepsilon
\end{equation}
provided that $k$ and $m$ are sufficiently large. By \eqref{409},
\eqref{412} and \eqref{416}, $\{ u_k\}$ is (up to subsequences) a
Cauchy sequence in measure in $\o$, and hence there exists a
measurable function $u:÷\o\rightarrow \R$ such that \eqref{408ter}
holds.
\medskip

\noindent{\bf Step 2.}
\begin{equation}\label{416bis}
\{\nabla u_k\} \hbox{\emph{ is a Cauchy sequence in measure}.}
\end{equation}
\par
Given any $t,\tau,\delta>0$, we have that
\begin{align}\label{417}
|\{|\nabla u_k-\nabla u_m|>t \}|& \le |\{|\nabla  u_k|>\tau
\}|+|\{| \nabla u_m|>\tau\}|+ |\{|u_k-u_m|>\delta  \}| \\
&+|\{|u_k-u_m|\le \delta,  |\nabla  u_k|\le \tau , |\nabla
u_m|\le \tau , |\nabla u_k-\nabla u_m|>t \}|,\nonumber
\end{align}
for $k,m\in \N$. Either assumption \eqref{401} or \eqref{402},
according to whether $q\in (1, \infty]$ or $q=1$, and Theorem
\ref{teo5ter} ensure, via \eqref{406}, that 
\begin{equation}\label{417bis}
\|\nabla u _k\|_{L^{p-1} (\o )} \leq C \|f\|_{L^q(\o
)}^{\frac{1}{p-1}}\,,
\end{equation}
for some constant $C$ independent of $k$. Hence,
\begin{equation}\label{418}
 |\{|\nabla  u_k|>\tau \}|\le \bigg( \frac{C\|f\|_{L^q(\o )}^{\frac 1{p-1}}}{\tau} \bigg)^{p-1},
\end{equation}
for $k\in \N$ and for some constant $C$ independent of $k$. Thus
$\tau$ can be chosen so large that
\begin{equation}\label{419}
 |\{|\nabla  u_k|>\tau \}|< \varepsilon,\qquad \hbox{for }k\in\N.
\end{equation}
Next, set
\begin{equation}\label{420}
G=\{ |u_k-u_m|\le \delta, |\nabla u_k|\le \tau,  |\nabla u_m|\le
\tau, |\nabla u_k -\nabla u_m|\ge t\}.
\end{equation}
We claim that, if \eqref{419} is fulfilled, then
\begin{equation}\label{421}
|G|< \varepsilon.
\end{equation}
To verify our claim, observe that, if we define
$$
S=\{(\xi, \eta)\in {\R} ^{2n}:÷|\xi|\le \tau, |\eta|\le \tau,
|\xi-\eta|\ge \tau \},
$$
and $l:÷\o\rightarrow [0,\infty)$ as
$$
l(x)=\inf\{[a(x, \xi)-a(x, \eta) ]\cdot (\xi-\eta): ÷(\xi,
\eta)\in S\},
$$
then $l(x)\ge 0$ and
\begin{equation}\label{422}
|\{l(x)=0\}|=0.
\end{equation}
Actually, this is a consequence of \eqref{4} and of the fact that
$S$ is compact and $a(x, \xi)$ is continuous in $\xi$ for every
$x$ outside a subset of $\o$ of Lebesgue measure zero.
\par\noindent
Now,
\begin{equation}\label{423}
\begin{split}
\int_G l(x)\, dx& \le\int_G [a(x, \nabla u_k)-a(x, \nabla u_m) ]\cdot (\nabla u_k-\nabla u_m)\, dx\\
&\le \int_{\{|u_k-u_m|\le \delta\}} [a(x, \nabla u_k)-a(x, \nabla u_m) ]\cdot (\nabla u_k-\nabla u_m)\, dx\\
&= \int_\o [a(x, \nabla u_k)-a(x, \nabla u_m) ]\cdot \nabla(T_{\delta}(u_k - u_m)) \, dx \\
&= \int_\o(f_k-f_m)T_{\delta}(u_k - u_m)dx\le 4
|\o|^\frac{1}{q'}\delta \|f\|_{L^q(\o)},
\end{split}
\end{equation}
where the last equality follows on making use of $T_{\delta}(u_k -
u_m)$ as test function in \eqref{408bis} for $k$ and $m$ and
substracting the resulting equations. Thanks to \eqref{422}, one
can show that for every $\varepsilon>0$ there exists $\theta>0$
such that if a measurable set $F\subset \o$ fulfills $\int_F
l(x)dx<\theta$, then $|F|<\varepsilon$. Thus, choosing $\delta$ so
small that $4 |\o|^\frac{1}{q'}\delta \|f\|_{L^q(\o)}<\theta$,
inequality \eqref{421} follows.

Finally,  since, by Step 1, $\{u_k\}$ is a Cauchy sequence in
measure in $\o$,

\begin{equation}\label{424}
|\{ |u_k-u_m|>\delta\}|<\varepsilon\,,
\end{equation}
if $k$ and $m$ are sufficiently large. Combining \eqref{417}, \eqref{419}, \eqref{421} and \eqref{424} yields
$$
|\{ |\nabla u_k-\nabla u_m|>t\}|<4\varepsilon\,,
$$
for sufficiently large $k$ and $m$. Property \eqref{416bis} is thus established.
\medskip

\noindent{\bf Step 3.} $u\in W^{1,p}_T(\o)$, \emph{and}
\begin{equation}\label{425}
\nabla u_k \rightarrow \nabla u \qquad \hbox{a.e. in }\o\,,
\end{equation}
\emph{up to subsequences, where $\nabla u$ is the generalized
gradient of $u$ in the sense of \eqref{504}. }

Since $\{\nabla u_k\}$ is a Cauchy sequence in measure, there exists a measurable function $Z:÷\o\rightarrow \rn$ such that
\begin{equation}\label{426}
\nabla u_k \rightarrow Z \qquad \hbox{a.e. in }\o
\end{equation}
(up subsequences). Fix any $t>0$. By \eqref{414}, $\{\T(u_k)\}$ is
bounded in $W^{1,p}(\o)$. Thus, there exists a function $\overline
u_t\in W^{1,p}(\o)$ such that
\begin{equation}\label{427}
\T(u_k)  \rightharpoonup  \overline u_t \qquad \hbox{weakly in
}W^{1,p}(\o)
\end{equation}
(up subsequences). By Step 1, $\T(u_k)\rightarrow \T(u)$ a.e. in $\o$, and hence
\begin{equation}\label{428}
  \overline u_t =\T (u) \qquad \hbox{a.e. in }\o.
  \end{equation}
Thus,
\begin{equation}\label{429}
\T(u_k)  \rightharpoonup \T (u) \qquad \hbox{weakly in
}W^{1,p}(\o).
\end{equation}
In particular, $u\in W^{1,p}_T(\o)$, and
\begin{equation}\label{430}
\nabla (\T(u))=\chi_{\{|u|<t\}}\nabla u  \qquad \hbox{a.e. in
}\o\,.
\end{equation}
By \eqref{408ter} and \eqref{426},
$$
\nabla (\T(u_k))=\chi_{\{|u_k|<t\}} \nabla u_k \rightarrow
\chi_{\{|u|<t\}}Z \qquad \hbox{a.e. in }\o\,.
$$
Hence, by \eqref{429}
\begin{equation}\label{431}
\nabla (\T(u))=\chi_{\{|u|<t\}}Z  \qquad \hbox{a.e. in }\o\,.
\end{equation}
Owing to the arbitrariness of $t$, coupling \eqref{430} and
\eqref{431} yields
\begin{equation}\label{432}
 Z =\nabla u\qquad \hbox{a.e. in }\o.
\end{equation}
Equation \eqref{425} is a consequence of \eqref{426} and \eqref{432}.
\smallskip
\par\noindent
{\bf Step 4.} $u\in V^{1,p-1}(\o)$ \emph{and satisfies property
$(i)$ of the definition of approximable solution}.

>From \eqref{425} and \eqref{417bis}, via Fatou's lemma, we deduce
that
$$
\|\nabla u\|_{L^{p-1}(\o)}\le C \|f\|_{L^q(\o)}^{\frac{1}{p-1}},
$$
for some constant $C$ independent of $f$.
Hence, $u\in V^{1,p-1}(\o)$.
\par\noindent
As far as property $(i)$ of the definition of approximable
solution is concerned, by \eqref{425}
\begin{equation}\label{433}
a(x, \nabla u_k) \rightarrow a(x, \nabla u)\qquad \hbox{for a.e.
$x \in \o$ .}
\end{equation}
Fix any $\Phi \in W^{1, \infty}(\o)$ and any measurable set
$F\subset \o$. Owing to Lemma \ref{lemnuovo} and \eqref{406},
\begin{align}\label{434}
\int_F |a(x, \nabla u_k)\cdot \nabla \phi|\, dx&\le \|\nabla
\phi\|_{L^\infty(\o)} \int_F \big(| \nabla u_k|^{p-1}+h(x)\big)\,
dx\\ \nonumber &\le\|\nabla \phi\|_{L^\infty(\o)} \left( \varsigma
(|F|) \|f\|_{L^q(\Omega )} + \int _F h(x)\, dx\right)
%
\end{align}
for some function $\varsigma : (0, \infty ) \to [0, \infty )$ such
that $\lim _{s\to 0^+}\varsigma (s)=0$. From \eqref{433} and
\eqref{434}, via Vitali's convergence theorem, we deduce that the
left-hand side of \eqref{408bis} converges to the left-hand side
of \eqref{231} as $k\rightarrow\infty$. The right-hand side of
\eqref{408bis} trivially converges to the right-hand side of
\eqref{231}, by \eqref{404}. This completes the proof of the
present step, and hence also the proof
of the existence of an approximable solution to \eqref{1}.
\smallskip
\par
We are now concerned  with the uniqueness of the solution to
\eqref{1}. Assume that $u$ and $\overline u$ are approximable
solutions to problem \eqref{1}. Then there exist sequences
$\{f_k\}$ and $\{\overline f_k\}\subset L^{q}(\o)\cap (V^{1,p}(\o
))'$ having the following properties: $\int _\Omega f_k dx = \int
_\Omega \overline f_k dx =0$; $f_k\rightarrow f$ and $\overline
f_k\rightarrow f$ in $L^q(\o)$; the weak solutions $u_k$ to
problem \eqref{407} and the weak solutions
 $\overline u_k$ to problem \eqref{407} with $f_k$ replaced by $\overline f_k$,
fulfill $u_k \rightarrow u$ and $\overline u_k \rightarrow
\overline u$ a.e. in $\o$. Fix any $t>0$ and choose the test
function $\Phi= \T(u_k-\overline u_k)$ in \eqref{408bis}, and in
the same equation with $u_k$ and $f_k$ replaced by $\overline u_k$
and $\overline f_k$, respectively. Subtracting the resulting
equations yields
\begin{equation}\label{435}
\int_\o\chi_{\{|u_k-\overline u_k|\le t\}} [a(x, \nabla u_k)-a(x,
\nabla \overline u_k) ]\cdot (\nabla u_k-\nabla \overline u_k)\,
dx=\int_\o (f_k-\overline f_k)\T(u_k-\overline u_k)\, dx
\end{equation}
for $k\in \N$. Since $|\T(u_k-\overline u_k)|\le t$ in $\o$ and
$f_k-\overline f_k\rightarrow 0$ in $L^q(\o)$, the right-hand side
of \eqref{435} converges to $0$ as $k\rightarrow\infty$. On the
other hand, arguments analogous to those exploited above in the
proof of the existence tell us that $\nabla u_k \rightarrow \nabla
u$ and $\nabla \overline u_k \rightarrow \nabla \overline u$ a.e.
in $\o$, and hence, by \eqref{4} and Fatou's lemma,
$$
\int_{\{|u-\overline u|\le t\}} [a(x, \nabla u)-a(x, \nabla
\overline u) ]\cdot (\nabla u-\nabla \overline u)\, dx=0.
$$
Thus, owing to \eqref{4}, we have that $\nabla u=\nabla \overline
u$ a.e. in $\{|u-\overline u|\le t\}$ for every $t>0$, and hence
\begin{equation}\label{436}
\nabla u=\nabla \overline u \qquad \hbox{a.e. in }\o.
\end{equation}
When $p\ge 2$, equation \eqref{436} immediately entails that
$u-\overline u=c$ in $\o$ for some $c\in \R$. Indeed, since
$u,\overline u \in V^{1, p-1}(\o)$ and $p-1\ge 1$, $u$ and
$\overline u$ are Sobolev functions in this case.
\par\noindent
The case when $1<p<2$ is more delicate. Consider a family
$\{\o_\varepsilon \}_{\varepsilon>0}$ of smooth open sets invading
$\o$. A version of the Poincare inequality \cite{Ma2, Z} tells us
that a constant $C(\o_\varepsilon)$ exists such that
\begin{equation}\label{437}
\left(\int_{\o_\varepsilon }|v-{\rm med }_{\o_\varepsilon}
(v)|^{n'}\, dx\right)^{\frac{1}{n'}}\le C(\o_\varepsilon)
\int_{\o_\varepsilon}|\nabla v|\, dx,
\end{equation}
for every $v\in W^{1,1}(\o_\varepsilon)$. Fix any $t,\tau>0$. An
application of \eqref{437} with $v= T_{\tau}(u-\T(\overline u))$,
and the use of \eqref{436} entail that
\begin{equation}\label{438}
\left( \int_{\o_\varepsilon }|T_{\tau}(u-\T(\overline u))-{\rm med
}_{\o_\varepsilon}(T_{\tau}(u-\T(\overline u))) |^{n'}\, dx \right
)^{\frac{1}{n'}}
\end{equation}
\begin{equation*}
\le C(\o_\varepsilon)
\left(
\int_{\{t<|u|<t+\tau\} }|\nabla u|\, dx+
\int_{\{t-\tau<|u|<t\} }|\nabla u|\, dx
\right ).
\end{equation*}
We claim that, for each $\tau>0$, the right-hand side of
\eqref{438} converges to 0 as $t\rightarrow\infty$. To verify this
claim, choose the test fnction $\Phi=T_{\tau}(u_k-\T(u_k))$ in
\eqref{408bis} and exploit \eqref{2} to deduce that
\begin{equation}\label{439}
\int_{\{t<|u_k|<t+\tau\}}|\nabla u_k|^p\, dx\le
\int_{\{t<|u_k|<t+\tau\}}a(x, \nabla u_k)\cdot \nabla u_k\, dx\le
\tau \int_{\{|u_k|>t\}}|f_k|\, dx.
\end{equation}
On passing to the limit as $k\rightarrow\infty$ in \eqref{439}
one can easily deduce that
\begin{equation}\label{440}
\int_{\{t<|u|<t+\tau\}}|\nabla u|^p\, dx\le\tau
\int_{\{|u|>t\}}|f|\, dx.
\end{equation}
Hence, the first integral on the right-hand side of \eqref{438}
approaches 0 as $t\rightarrow \infty$. An analogous argument
shows
that also the last integral in \eqref{438} goes to 0 as
$t\rightarrow \infty$. Since
$$
\lim_{t \to \infty} \Big( T_{\tau}(u-\T(\overline u))-{\rm med
}_{\o_\varepsilon}(T_{\tau}(u-\T(\overline u)))\Big) =
T_{\tau}(u-\overline u)-{\rm med }_{\o_\varepsilon}
(T_{\tau}(u-\overline u)), \quad \hbox{a.e. in } \o,
$$
from \eqref{438}, via Fatou's lemma, we obtain that
\begin{equation}\label{441}
\int_{\o_\varepsilon}|T_{\tau}(u-\overline u)-{\rm med
}_{\o_\varepsilon}(T_{\tau}(u-\overline u)) |^{n'}\,dx=0
\end{equation}
for $\tau>0$. Thus, the integrand in \eqref{441} vanishes a.e. in
$\o_\varepsilon$ for every  $\tau >0$, and hence also  its limit
as $\tau \rightarrow \infty$ vanishes a.e. in $\o_\varepsilon$.
Therefore,
a constant $c(\varepsilon)$ exists such that $u-\overline
u=c(\varepsilon)$ in $\o_\varepsilon$ for every $\varepsilon >0$.
Consequently, $u-v=c$ in $\o$ for some $c\in\R$.

\section{Strongly monotone operators}\label{sec4}

\subsection{Continuous dependence estimates}\label{subsec4.1}
The present subsection is concerned with a norm estimate for the
difference of the gradients of weak solutions to problem
\eqref{1}, with different right-hand sides in $(V^{1,p}(\o ))'$,
under the strong monotonicity assumption \eqref{5}. Such an
estimate is a crucial ingredient for a variant (a simplification
in fact) in the approach to existence presented in Section
\ref{sec3}, and leads to the continuous dependence result of
Theorem \ref{teo4}.

\begin{theorem}\label{teo5}

Let $\o$, $p$, $q$, $a$, $f$ and $g$ be as in Theorem \ref{teo4}.
Assume, in addition, that $f,g\in  (V^{1,p}(\o ))'$. Let $u$ be a
weak solution to problem \eqref{1}, and let $v$ be a weak solution
to problem \eqref{1} with $f$ replaced by $g$. Let $0<\sigma\le p$
and let $r=\max \{p,2 \}$. Then, there exists a constant $C$ such
that
\begin{equation}\label{501}
\|\nabla u-\nabla v \|_{L^\sigma (\o )} \leq C \|f-g\|_{L^q(\o
)}^{\frac{1}r} \left(\|f\|_{L^q(\o)}+\|g\|_{L^q(\o)}
\right)^{\frac{1}{p-1}-\frac{1}{r}}, \,
\end{equation}
if either
\par\noindent
(i) $q>1$, $q(p-1) \leq \sigma $ and \eqref{5teri} holds,
\par\noindent
or
\par\noindent
(ii) $1<q<\infty$, $0 < \sigma < q(p-1)$ and \eqref{5terii} holds,
\par\noindent
or
\par\noindent
(iii)  $q =\infty$ and \eqref{5teriii} holds,
\par\noindent
or
\par\noindent
(iv)  $q =1$ and \eqref{5teriv} holds.
\par\noindent
Moreover, the constant $C$ in \eqref{501} depends only on $p$,
$q$, $\sigma$ and on the left-hand side either of \eqref{5teri},
or \eqref{5terii}, or \eqref{5teriii}, or \eqref{5teriv},
respectively.
\end{theorem}
\medskip

\noindent{\bf Proof.} Throughout the proof, $C$ and $C'$ will
denote constants which may change from equation to equation, but
which depend only on the quantities specified in the statement.
\par\noindent  We shall focus  on the case
where $\sigma<p,$ the case where $\sigma=p$ being analogous, and
even simpler. \par First, assume that $1<q<\infty$, and hence that
we are dealing either with case $(i)$ or $(ii).$ We may suppose,
without loss of generality, that $u$ and $v$ are normalized in
such a way that ${\rm med }_\o (u)={\rm med }_\o (v)=0$. Given any
$\gamma \in (-1 , 0)$ and $\varepsilon > 0$, define
$(u-v)_{\varepsilon , \gamma } : \o \to \R$ as
$$
(u-v)_{\varepsilon , \gamma }=\max\{(u-v)_+\,, \varepsilon
\}^{\gamma + 1} - \max\{(u-v)_-\,, \varepsilon  \}^{\gamma + 1}.
$$
%
%
%
The chain rule for derivatives in Sobolev spaces ensures that
$(u-v)_{\varepsilon , \gamma }\in V^{1,p}(\o)$, and that
$$
\nabla (u-v)_{\varepsilon , \gamma }= (\gamma +1 )
|u-v|^{\gamma}\chi_{|u-v|>\varepsilon}\nabla (u-v) \qquad \quad
\hbox{a.e. in $\o$.}
$$
Thus,
the function $(u-v)_{\varepsilon, \gamma }$ can be used as test
function $\Phi$ is the definition of weak solution for $u$ and
$v$. Subtracting the resulting equations yields
\begin{equation}\label{508}
(\gamma+1)\int_{\{|u-v|>\varepsilon  \}}|u-v|^\gamma [a(x, \nabla
u)- a(x, \nabla v)]\cdot \nabla (u-v)\, dx=
\int_\o(u-v)_{\varepsilon, \gamma } (f-g)\, dx.
\end{equation}
If $1<p<2$,  making use of \eqref{5} and passing to the limit as
$\varepsilon \rightarrow 0^+$ in \eqref{508} tell us that
\begin{equation}\label{509}
\int_\o |u-v|^{\gamma} \frac{|\nabla u-\nabla v|^2}{(|\nabla
u|+|\nabla v|)^{2-p}}\, dx\le C \int_\o |u-v|^{\gamma+1} |f-g|\,
dx\,.
\end{equation}
If $p\ge 2$, then the same argument
yields
\begin{equation}\label{510}
\int_\o |u-v|^{\gamma}|\nabla u-\nabla v|^p\, dx \le C\int_\o
|u-v|^{\gamma+1} |f-g|\, dx.
\end{equation}
\par
Consider first the case when $1<p<2$. Let $(\alpha, \beta, \varrho
)$ be the solution to the system
\begin{equation}\label{511}
(p-2)\frac{\beta}{2-\beta}=\sigma,
\end{equation}
\begin{equation}\label{512}
\frac{1-\alpha}{\alpha}\frac{\sigma\beta}{\beta-\sigma}=\varrho,
\end{equation}
\begin{equation}\label{513}
\frac{\alpha\varrho}{\alpha\varrho-2\alpha+1}=q.
\end{equation}
Namely,
\begin{equation}\label{515}
\alpha=\frac{\sigma q +pq-2\sigma}{2pq-2\sigma}\,,
\end{equation}
\begin{equation}\label{514}
\beta=\frac{2\sigma}{2+q-p}\,,
\end{equation}
\begin{equation}\label{516}
\varrho=\frac{2\sigma q}{q\sigma+pq-2\sigma}.
\end{equation}
Observe that $\alpha \in (1/2, 1)$, $\beta>\sigma$ and
\begin{equation}\label{517}
\alpha\varrho=\frac{\sigma q }{pq-\sigma}.
\end{equation}
Now, choose $\gamma=2(\alpha-1)$ in \eqref{509}, and note that
actually $\gamma \in (-1, 0)$. Thus, the following chain holds:
\begin{align}\label{518}
\int_\o&\left(|\nabla u-\nabla v||u-v| ^{\alpha-1}\right)^\beta \,
dx   \le \left(\int_\o\frac{(|\nabla u-\nabla v||u-v|
^{\alpha-1})^2}{(|\nabla u|+|\nabla v|)^{2-p}} \, dx
\right)^{\frac{\beta}{2}}  \left(\int_\o(|\nabla u|+|\nabla
v|)^{\sigma} \, dx \right)^{1-\frac{\beta}{2}} \intertext{\qquad
\qquad \qquad \qquad \qquad \qquad \qquad \qquad \qquad \qquad
\qquad \qquad \qquad \qquad \qquad (by H\"older's inequality)}
\nonumber
 \nonumber
&\le  \left(\int_\o |f-g||u-v| ^{2(\alpha-1)+1} \, dx
\right)^{\frac{\beta}{2}} \left(\int_\o(|\nabla u|+|\nabla
v|)^{\sigma} \, dx \right)^{1-\frac{\beta}{2}} \intertext{\qquad
\qquad \qquad \qquad \qquad \qquad \qquad \qquad \qquad \qquad
\qquad \qquad \qquad \qquad \qquad (by \eqref{509} and
\eqref{511})} \nonumber
 \nonumber
&\le
 \left(\int_\o |u-v| ^{\alpha\varrho}
\, dx \right)^{\frac{2(\alpha-1)+1}{\alpha\varrho}\frac{\beta}{2}}
\left(\int_\o |f-g|
^{\frac{\alpha\varrho}{\alpha\varrho-2(\alpha-1)-1}} \, dx
\right)^{\frac{\alpha \varrho-2(\alpha-1)-1}{\alpha
\varrho}\frac{\beta}{2}} \left(\int_\o(|\nabla u|+|\nabla
v|)^{\sigma} \, dx \right)^{1-\frac{\beta}{2}} \intertext{\qquad
\qquad \qquad \qquad \qquad \qquad \qquad \qquad \qquad \qquad
\qquad \qquad \qquad \qquad \qquad (by H\"older's inequality)}
\nonumber &= \left(\int_\o |u-v| ^{\frac{\sigma q}{pq-\sigma}} \,
dx \right)^{\frac{\beta}{2q'}} \left(\int_\o |f-g| ^{q} \, dx
\right)^{\frac{\beta}{2q}} \left(\int_\o(|\nabla u|+|\nabla
v|)^{\sigma} \, dx \right)^{1-\frac{\beta}{2}}\nonumber
\intertext{\qquad \qquad \qquad \qquad \qquad \qquad \qquad \qquad
\qquad \qquad \qquad \qquad \qquad \qquad \qquad (by \eqref{513}
and \eqref{517}).} \nonumber
\end{align}
Next, observe that $\sigma \ge q(p-1)$ if and only if
$\frac{\sigma q}{pq-\sigma}\ge q(p-1)$. Thus, an application of
Theorem \ref{teoC} with $\sigma$ replaced by $\frac{\sigma
q}{pq-\sigma}$ tells us that, either under \eqref{502} or
\eqref{503}, according to whether $\sigma \ge q(p-1)$ or $\sigma <
q(p-1)$, one has that
\begin{equation}\label{519}
\|u\|_{L^{\frac{\sigma q}{pq-\sigma}}(\o)}\le C
\|f\|_{L^q(\o)}^{\frac{1}{p-1}} \qquad \hbox{and} \qquad
\|v\|_{L^{\frac{\sigma q}{pq-\sigma}}(\o)}\le C
\|g\|_{L^q(\o)}^{\frac{1}{p-1}}\,.
\end{equation}
Moreover, by Theorem \ref{teo5ter},
\begin{equation}\label{520}
\|\nabla u\|_{L^{\sigma}(\o)}\le C \|f\|_{L^q(\o)}^{\frac{1}{p-1}}
\qquad \hbox{and} \qquad \|\nabla v\|_{L^{\sigma}(\o)}\le C
\|g\|_{L^q(\o)}^{\frac{1}{p-1}}\,.
\end{equation}
Combining \eqref{518}-\eqref{520} yields
\begin{equation}\label{521}
\int_\o|\nabla |u-v|^\alpha|^\beta\, dx\le C \left(\int_\o |f-g|
^{q} \, dx \right)^{\frac{\beta}{2q}}
\left(\|f\|_{L^q(\o)}+\|g\|_{L^q(\o)}
\right)^{\frac{\sigma}{p-1}(1-\frac{\beta}{2})+\frac{\sigma(q-1)}{(p-1)(pq-\sigma)}\frac{\beta}{2}}.
\end{equation}
By Holder's inequality, \eqref{512},  \eqref{517} and  \eqref{519},
\begin{align}\label{522}
\int_\o |\nabla u-\nabla v| ^{\sigma} \, dx   &= \frac{1}{\alpha^\sigma}\int_\o|\nabla |u-v|^\alpha|^\sigma |u-v|^{(1-\alpha)\sigma}\, dx \\
&\le \frac{1}{\alpha^\sigma}\left(\int_\o|\nabla
|u-v|^\alpha|^\beta\, dx   \right)^\frac{\sigma}{\beta}
\left(\int_\o |u-v| ^{\frac{\sigma q}{pq-\sigma}} \, dx
\right)^{1-\frac{\sigma}{\beta}}\nonumber
 \\
&\le C\left(\int_\o|\nabla |u-v|^\alpha|^\beta\, dx
\right)^\frac{\sigma}{\beta} \left(\|f\|_{L^q(\o)}+\|g\|_{L^q(\o)}
\right)^{\frac{(\beta-\sigma)\sigma
q}{(p-1)\beta(pq-\sigma)}}\,.\nonumber
\end{align}
Inequality \eqref{501} follows from \eqref{521} and \eqref{522}.

Assume now that $p\ge 2$. Let $(\alpha , \varrho )$ be the
solution to the system
\begin{equation}\label{523}
\frac{\alpha\varrho}{\alpha\varrho-p(\alpha-1)-1} =q,
\end{equation}
\begin{equation}\label{524}
\frac{1-\alpha }{\alpha }\frac{\sigma p}{p-\sigma}=\varrho,
\end{equation}
namely
$$
\varrho=\frac{q\sigma p}{pq(p-1)+\sigma(q-p)} \qquad
\quad\hbox{and}\qquad \quad \alpha
=\frac{pq(p-1)+\sigma(q-p)}{p(pq-\sigma)}.
$$
In particular,
\begin{equation}\label{525}
\alpha \varrho=\frac{\sigma q}{pq-\sigma}
\end{equation}
also in this case. Take $\gamma = p(\alpha -1)$ in \eqref{510}, an
admissible choice since $p(\alpha -1) \in (-1, 0)$. From
\eqref{510}, \eqref{525}, \eqref{523} and \eqref{519} one deduces
that
\begin{align}\label{526}
\int_\o |\nabla |u-v|^\alpha|^p \, dx &= \alpha ^p\int_\o
\left(|\nabla u-\nabla v||u-v| ^{\alpha-1}\right)^p\, dx
\\ \nonumber
&\le C\int_\o |f-g||u-v| ^{p(\alpha-1)+1}
\, dx
\\ \nonumber
&\le C\left(\int_\o |u-v| ^{\frac{\sigma q}{pq-\sigma}} \, dx
\right)^{\frac{1}{q'}} \left(\int_\o |f-g| ^{q} \, dx
\right)^{\frac{1}{q}}
 \\ \nonumber
&\le C' \left(\|f\|_{L^q(\o)}+\|g\|_{L^q(\o)}
\right)^{\frac{\sigma(q-1)}{(p-1)(pq-\sigma)}}\|f-g\|_{L^q(\o)}\,.
\end{align}
Analogously to \eqref{522}, we have
that
\begin{align}\label{527}
\int_\o |\nabla u-\nabla v| ^{\sigma} \, dx   &\le
  \frac{1}{\alpha^\sigma}\left(\int_\o|\nabla |u-v|^\alpha|^\beta\, dx   \right)^\frac{\sigma}{p}
\left(\int_\o |u-v| ^{\frac{\sigma q}{pq-\sigma}} \, dx
\right)^{1-\frac{\sigma}{p}}
 \\ \nonumber
&\le C\left(\int_\o|\nabla |u-v|^\alpha|^p\, dx
\right)^\frac{\sigma}{p} \left(\|f\|_{L^q(\o)}+\|g\|_{L^q(\o)}
\right)^{\frac{\sigma q}{(p-1)(pq-\sigma)}(1-\frac{\sigma}{p})}\\
\nonumber &\le C' \left(\|f\|_{L^q(\o)}+\|g\|_{L^q(\o)}
\right)^{\frac{\sigma}{p(p-1)}}\|f-g\|_{L^q(\o)}^\frac{\sigma}{p},
\end{align}
where the second inequality holds
owing to \eqref{519} and the last one  to \eqref{526}. This
completes the proof of \eqref{501} in cases $(i)$ and $(ii)$.

Case $(iii)$ can be dealt with an analogous argument, requiring
easy modifications. The details are omitted for brevity.

Finally, consider case $(iv)$.  As above, we may assume that ${\rm
med }_\o (u)={\rm med }_\o (v)=0$. Let us set
$$w=(u-v)_+.$$
Given any integrable function $\zeta : (0,\m2) \rightarrow
[0,\infty)$, define  $\Lambda:÷[0, \m2]\rightarrow[0,\infty)$ as
\begin{equation}\label{541}
\Lambda(r)=\int_0^r \zeta (\rho)\, d\rho, \qquad \hbox{for }r\in
[0, \m2].
\end{equation}
Moreover, for any fixed $s\in [0, \m2]$, define $I:÷[0,
\m2]\rightarrow [0,\infty)$ as
\begin{equation}\label{542}
I(r) =\left\{
\begin{array}{ccc}
\Lambda(r) & \text{if} &0\le r\le s, \\
\Lambda(s)& \text{if} &s<r\leq \m2,
\end{array}
\right.
\end{equation}
and $\Phi:÷\o\rightarrow [0, \infty )$ as
\begin{equation}\label{543}
\Phi(x)=\int_0^{w(x)}I (\mu_{w}(t))\, dt, \qquad \hbox{for
}x\in\o.
\end{equation}
Since $I\circ \mu _{w}$ is a bounded function, the chain rule for
derivatives
 in Sobolev spaces
tells us that $\Phi\in V^{1,p}(\o)$ and
\begin{equation}\label{544}
\nabla \Phi=\chi_{\{u-v>0\}}I (\mu_{w}(w))(\nabla u-\nabla
v)\qquad \hbox{a.e. in $\o$}.
\end{equation}
Choosing $\Phi$ as test function in the definitions of weak
solution for $u$ and $v$ and subtracting the resulting equations
yields
\begin{equation}\label{545}
\int_{\{u-v>0\}}I(\mu_{w}(w(x)))[a(x, \nabla u)-a(x, \nabla v)
]\cdot (\nabla u-\nabla v)\, dx=\int_\o (f-g)\Phi\, dx.
\end{equation}
Observe that
\begin{align}\label{546}
\|\Phi\|_{L^\infty(\o)}&\le \int_0^\infty I (\mu_{w}(t))\, dt\\
\nonumber  &=\int^\infty_{w^*(s)}\Lambda(\mu_{w}(t))\,
dt+\int_0^{w^*(s)}\Lambda(s)\, dt
\\ \nonumber
&=\int^\infty_{w^*(s)}\int_0^{\mu_{w}(t)}\zeta(\rho)d\rho\,
dt+\Lambda(s)w^*(s)
\\ \nonumber
&=\int_0^s(w^*(\rho)-w^*(s))\zeta(\rho)\, d\rho+\Lambda(s)w^*(s)
\\ \nonumber
&=\int_0^s w^*(\rho)\zeta(\rho)\, d\rho
\\ \nonumber
&\le \int_0^s (u_+^*(\rho/2)+v_-^*(\rho/2))\zeta(\rho)\, d\rho
\\ \nonumber
&\le\left(\|f_+\|_{L^1(\o)}^\frac{1}{p-1}+\|g_-\|_{L^1(\o)}^\frac{1}{p-1}
\right)\int_0^s\nu _p(\rho /2)^\frac{1}{1-p}\zeta(\rho)\, d\rho,
\end{align}
where the third equality holds by Fubini's theorem, the last but
one inequality by \eqref{somma}, and the last inequality by
estimate \eqref{A} and by the corresponding estimate for $v$.
Combining \eqref{545} and \eqref{546} entails that
\begin{align}\label{547}
\int_{\{u-v>0\}} & I(\mu_{w}(w(x))) [a(x, \nabla u)-a(x, \nabla v)
]\cdot (\nabla u-\nabla v)\, dx \\ \nonumber & \le
\|(f-g)_+\|_{L^1(\o)}\left(\|f_+\|_{L^1(\o)}^\frac{1}{p-1}+\|g_-\|_{L^1(\o)}^\frac{1}{p-1}
\right)\int_0^s \nu _p(r/2)^\frac{1}{1-p}\zeta(r)\, dr.
\end{align}
Let us distinguish the cases when $1<p<2$ and $p \geq 2$.
\par\noindent
First, assume that $p\ge 2$. By \eqref{5} and \eqref{547},
\begin{multline}\label{548}
C\int_{\{u-v>0\}}|\nabla w|^pI (\mu_{w}(w(x)))\, dx
\\
\le
\|(f-g)_+\|_{L^1(\o)}\left(\|f_+\|_{L^1(\o)}^\frac{1}{p-1}+\|g_-\|_{L^1(\o)}^\frac{1}{p-1}
\right)\int_0^s \nu _p(r/2)^\frac{1}{1-p}\zeta(r)\, dr\,.
\end{multline}
Since $w$ and $w^*$ are equidistributed functions and $I$ is
non-decreasing,
\begin{align}\label{549}
\big(I \circ \mu_{w} \circ w \big)_*(r)= \big(I \circ \mu_{w}
\circ w^* \big)_*(r) \ge I _*(r) =I (r), \quad \hbox{for }r\in(0,
|\o |/2)\,.
\end{align}
Hence, by \eqref{HL},
\begin{align}\label{550}
\int_{\{u-v>0\}}|\nabla w|^p I (\mu_{w}(w(x)))\, dx
&\ge\int_0^{\M2}|\nabla w|^*(r)^p\big(I \circ \mu_{w} \circ w^*
\big)_*(r)\, dr
\\
&\ge \int_0^{\M2}|\nabla w|^*(r)^pI (r)\, dr\nonumber
\\
&\ge \int_0^s|\nabla w|^*(r)^pI (r)\, dr\nonumber
\\
&\ge |\nabla w|^*(s)^p\int_0^s\int_0^r\zeta(\rho)\,\ d\rho\,
dr\quad \hbox{for }s\in(0,\m2).\nonumber
\end{align}
>From \eqref{548} and \eqref{550} we obtain that
\begin{equation}\label{551}
C |\nabla w|^*(s)^p\frac{\int_0^s\zeta(r)(s-r)\,\
dr}{\int_0^s\nu _p(r/2)^\frac{1}{1-p}\zeta(r)\, dr }\le
\|(f-g)_+\|_{L^1(\o)}\left(\|f_+\|_{L^1(\o)}^\frac{1}{p-1}+\|g_-\|_{L^1(\o)}^\frac{1}{p-1}
\right),
\end{equation}
for $s\in(0,\m2)$. Clearly,
\begin{equation}\label{552}
\sup_{\zeta}\frac{\int_0^s\zeta(r)(s-r)\,\ dr}{\int_0^s\nu
_p(r/2)^\frac{1}{1-p}\zeta(r)\, dr }=\|(s-r)\nu
_p(r/2)^\frac{1}{p-1}\|_{L^\infty(0,s)}\ge \frac s2 \nu
_p(s/4)^\frac{1}{p-1},
\end{equation}
for $s\in(0,\m2)$. Thus, owing to the arbitrariness of $\zeta$,
inequality \eqref{551} implies that
\begin{equation}\label{553}
 |\nabla w|^*(s)^\sigma\le \frac{C}{\big( s \nu _p(s/4)^\frac{1}{p-1}
  \big)^{\frac{\sigma}{p}}}\left( \|(f-g)_+\|_{L^1(\o)}\left(\|f_+\|_{L^1(\o)}^\frac{1}{p-1}
  +\|g_-\|_{L^1(\o)}^\frac{1}{p-1} \right) \right)^{\frac{\sigma}{p}},
\end{equation}
for $s\in(0,\m2)$. An analogous argument yields a similar
inequality with $(u-v)_+$ (i.e. $w$) replaced by $(u-v)_-$. Hence,
by \eqref{5teriv}, inequality \eqref{501} follows.
\par
Consider now the case when $1<p<2$. Define $H:÷\o\rightarrow
[0,\infty) $ by
$$
H(x)=\frac{|\nabla w|^{\frac{2}{p}}}{(|\nabla u|+|\nabla
v|)^{\frac{2-p}{p}}}, \quad \hbox{for }x\in\o.
$$
>From \eqref{5} and \eqref{547} we deduce that
\begin{align}\label{554}
C\int_{\{u-v>0\}}&H (x) I (\mu_{w}(w(x)))\, dx \\ &\le
 \|(f-g)_+\|_{L^1(\o)}\left(\|f_+\|_{L^1(\o)}^\frac{1}{p-1}+\|g_-\|_{L^1(\o)}^\frac{1}{p-1}
 \right)\int_0^s\nu _p(r/2)^\frac{1}{1-p}\zeta(r)\,
 dr\,.\nonumber
\end{align}
The same
argument leading to \eqref{553} now shows that
\begin{equation}\label{555}
H^* (s)\le \frac{C}{\big( s \nu _p(s/4)^\frac{1}{p-1} \big
)^{\frac{1}{p}}}\left(
\|(f-g)_+\|_{L^1(\o)}\left(\|f_+\|_{L^1(\o)}^\frac{1}{p-1}+\|g_-\|_{L^1(\o)}^\frac{1}{p-1}
\right) \right)^{\frac{1}{p}}
\end{equation}
for $s\in (0, |\o|/2)$. On the other hand,
\begin{align}\label{556}
|\nabla w|^*(s) &=\left (H^\frac{p}{2}(|\nabla u|+|\nabla
v|)^{\frac{2-p}{2}}\right )^*(s)
\\
&\le H^*(s/2)^\frac{p}{2}\left(|\nabla u|+|\nabla
v|\right)^*(s/2)^{\frac{2-p}{2}}\nonumber
\\
&\le CH^*(s/2)^\frac{p}{2} \left(|\nabla
u|^*(s/4)^{\frac{2-p}{2}}+|\nabla
v|^*(s/4)^{\frac{2-p}{2}}\right)\nonumber
\\
&\le C' H^*(s/2)^\frac{p}{2} \left( s \nu _p(s/8)^\frac{1}{p-1}
\right )^{\frac{p-2}{2p}}
\left(\|f\|_{L^1(\o)}^\frac{1}{p-1}+\|g\|_{L^1(\o)}^\frac{1}{p-1}
\right)^{\frac{2-p}{2}}\nonumber
\end{align}
for $s\in (0, |\o|/2)$.
Note that the first inequality holds by \eqref{prodotto}, the
second one by \eqref{somma} and last one by Proposition
\ref{prop7}. Coupling \eqref{555} and \eqref{556} yields
\begin{equation}\label{557}
 |\nabla w|^*(s)^\sigma\le \frac{C}{\left( s \nu _p(s/8)^\frac{1}{p-1} \right )^{\frac{\sigma}{p}}}
 \|(f-g)_+\|_{L^1(\o)}^{\frac{\sigma}{2}}\left(\|f\|_{L^1(\o)}+\|g\|_{L^1(\o)} \right)^{\frac{(3-p)\sigma}{2(p-1)}}
\end{equation}
for $s\in (0, \m2)$.
Inequality \eqref{557}, and a similar inequality with $(u-v)_+$
(i.e. $w$) replaced by $(u-v)_-$, imply \eqref{501} when
\eqref{5teriv} is in force.

\subsection{Proof of Theorem \ref{teo4}}\label{subsec4.2}

We proceed through the same steps and make use of the same
notations as in the proof of Theorem \ref{teo4bis}. The proofs of
Steps 1 and 3 are exactly the same. Thus, we shall focus on Steps
2 and 4.

\noindent{\bf Step 2.}  $\{\nabla u_k \}$ \emph{is a Cauchy
sequence in measure}.
\smallskip
\par
Either assumption \eqref{401} or \eqref{402}, according to whether
$q\in (1, \infty]$ or $q=1$, Theorem \ref{teo5} and \eqref{406}
ensure that
\begin{equation}\label{450}
\|\nabla u_k-\nabla u_m\|_{L^{p-1}(\o)}\le C
\|f_k-f_m\|_{L^{q}(\o)}^{\frac{1}{r}}\|f
\|_{L^{q}(\o)}^{\frac{1}{p-1}-\frac{1}{r}}
\end{equation}
for $k,m \in \N$ and for some constant $C$ independent of $k$ and
$m$, where $r=\max\{ p,2\}$.
Hence, $\{\nabla u_k \}$ is a Cauchy sequence in measure.

\noindent{\bf Step 4.} $u\in V^{1, p-1}(\o)$, \emph{and satisfies
property $(i)$ of the definition of approximable solution}.
\smallskip
\par
>From Theorem \ref{teo5ter}, Step 3 and Fatou's Lemma, we get that
$$
\|\nabla u \|_{L^{p-1}(\o)}\le \liminf_{k\to \infty}\|\nabla u_k
\|_{L^{p-1}(\o)}\le C\|f\|_{L^{q}(\o)},
$$
for some constant $C$, whence  $u\in V^{1, p-1}(\o)$. In order to prove \eqref{231}, note that
\begin{equation}\label{452}
\left|\frac{\partial u_k}{\partial x_i}\left|\frac{\partial
u_k}{\partial x_i}\right|^{p-2} - \frac{\partial u_m}{\partial
x_i}\left|\frac{\partial u_m}{\partial
x_i}\right|^{p-2}\right|\leq 2^{2-p} \left|\frac{\partial
u_k}{\partial x_i}-\frac{\partial u_m}{\partial
x_i}\right|^{p-1}\,,\qquad i=1,...,n,
\end{equation}
for $k,m\in \N$. Coupling \eqref{450} and \eqref{452} entails that
$ \left\{ \frac{\partial u_k}{\partial x_i}\left|\frac{\partial
u_k}{\partial x_i}\right|^{p-2}\right\}$ is a Cauchy sequence in
$L^{1}(\o)$, \hbox{for }i=1,...,n. Thus, the sequence  $\left\{
\frac{\partial u_k}{\partial x_i}\left|\frac{\partial
u_k}{\partial x_i}\right|^{p-2}\right\}$ converges to some
function in $L^{1}(\o)$. Since $\nabla u \rightarrow \nabla u $
a.e. in $\o$ by Step 3, necessarily
\begin{equation}\label{453}
\frac{\partial u_k}{\partial x_i}\left|\frac{\partial
u_k}{\partial x_i}\right|^{p-2}  {\rightarrow} \frac{\partial
u}{\partial x_i}\left|\frac{\partial u}{\partial
x_i}\right|^{p-2}\, \quad \hbox{in }L^{1}(\o), \quad \hbox{for
}i=1,...,n.
\end{equation}
Now, define the Carath\'eodory function $b:÷\o\times
\R^n\rightarrow\R^n$ as
\begin{equation}\label{454}
b(x, \eta )=a(x, \eta_1|\eta_1|^{\frac{2-p}{p-1}}, \dots ,
\eta_n|\eta_n|^{\frac{2-p}{p-1}}), \qquad \hbox{for }(x, \eta)\in
\o\times \rn .
\end{equation}
Hence,
\begin{equation}\label{455}
a(x,\xi)=b(\xi_1|\xi_1|^{p-2}, \dots ,\xi_n|\xi_n|^{p-2}), \qquad
\hbox{for }(x, \xi)\in \o\times \rn ,
\end{equation}
and, by \eqref{3}, for a.e. $x \in \o$
\begin{equation}\label{456}
|b(x, \eta)|\le C(|\eta|+h(x)), \qquad \hbox{for } \eta \in \rn\,.
\end{equation}
By \eqref{456}, the Nemytski operator $N:
÷(L^1(\o))^n\rightarrow(L^1(\o))^n$, defined by $Nz(x)=b(x,z(x))$,
for $z\in (L^1(\o))^n$, is continuous (see \cite[Section
26.3]{Zeidler}). Thus, by \eqref{453} and \eqref{455},
$$
 a(x, \nabla u_k) \rightarrow a(x, \nabla u)\qquad \hbox{in }
 (L^1(\o))^n\,.
 $$
 Consequently,
 $$
\lim_{k\to \infty} \int_\Omega a(x, \nabla u_k) \cdot \nabla \Phi
\, dx=\int_\Omega   a(x, \nabla u) \cdot \nabla \Phi \, dx
 $$
 for every $ \Phi \in W^{1,\infty}(\o)$. Trivially,
 $$
 \lim_{k\to \infty}\int_\Omega f_k \Phi \, dx = \int_\Omega f \Phi \,
 dx\,,
 $$
 for any such $ \Phi$. Hence, \eqref{231} follows from \eqref{408bis}. This completes
  the proof of the existence of an approximable solution to \eqref{1}.

As far as \eqref{403} is concerned, by the definition of
approximable solution, there exists sequences $\{f_k\}$ and
$\{g_k\}$ in $L^q(\o)\cap (V^{1,p}(\o ))'$ such that $f_k \to f$
and $g_k \to g$ in $L^q(\o )$ as $k \to \infty$, and such that the
sequences $\{u_k\}$ and $\{v_k\}$  of the weak solutions to
problem \eqref{1}, with $f$ replaced by $f_k$ and $f$ replaced by
$g_k$, converge to $u$ and $v$, respectively, a.e. in $\o$.
%
>From  Theorem \ref{teo5} we have that
\begin{equation}\label{457}
\|\nabla u_k-\nabla v_k\|_{L^{p-1}(\o)}\le
C\|f_k-g_k\|_{L^{q}(\o)}^{\frac{1}{r}} \left(\|f_k\|_{L^q(\o)}
+\|g_k\|_{L^q(\o)}   \right)^{\frac{1}{p-1}-\frac{1}{r}}
\end{equation}
for $k\in \N$. The same argument as in the proof of existence
above tells us that $\nabla u_k \rightarrow \nabla u $ and $\nabla
v_k \rightarrow \nabla v $ a.e. in $\o$ (up to subsequences).
Hence, by Fatou's lemma, we deduce \eqref{403}.

In particular, if $f=g$, then \eqref{403} entails that $\nabla u=
\nabla v$ a.e. in $\o$. Thus, the same argument as in the proof of
Theorem \ref{teo4bis} ensures that $u-v=c$ in $\o$ for some $c\in
\R$. This establishes the uniqueness of the solution up to
additive constants. \qed
\section{Applications and examples}\label{sec5}
Before presenting some applications of Theorems \ref{teo4bis} and
\ref{teo4} to special domains and classes of domains $\o $,
we state, for comparison, counterparts of these results involving
 the isoperimetric function $\lambda$. They immediately follow
from Theorems \ref{teo4bis} and \ref{teo4}, via \eqref{-514}.

\begin{corollary}\label{cor3bis}
Let $\o$, $p$, $q$, $a$ and $f$ be as in Theorem \ref{teo4bis}.
Assume that either
\par\noindent {\rm (i)} $1<q \leq \infty $ and
\begin{equation}\label{+301}
\int _0^{\m2}s^{\frac {q'}p}\bigg(\int _s^{\m2} \frac{dr}{\lambda
(r)^{p'}}\bigg)^{\frac {q'}{p'}} ds < \infty\,,
\end{equation}
or
\par\noindent
{\rm (ii)} $q=1$ and
\begin{equation}\label{+302}
\int _0^{\m2}s^{\frac 1p}\bigg(\int _s^{\m2} \frac{dr}{\lambda
(r)^{p'}}\bigg)^{\frac {1}{p'}} \frac{ds}{s} < \infty\,.
\end{equation}
Then there exists a unique (up to additive constants) approximable
solution to problem \eqref{1}.
\end{corollary}

\begin{corollary}\label{cor3}
Let $\o$, $p$, $q$, $a$ and $f$ be as in Theorem \ref{teo4}.
Assume that either $1<q \leq \infty $ and \eqref{+301} holds, or
$q=1$ and \eqref{+302} holds.
%
%
\par\noindent
Then there exists a unique (up to additive constants) approximable
solution to problem \eqref{1}  depending continuously on the
right-hand side of the equation. Precisely, if $g$ is another
function from $L^q(\o )$ such that $\int _\o g(x) dx =0$, and $v$
is the solution to \eqref{1} with $f$ replaced by $g$, then
\begin{equation*}
\|\nabla u - \nabla v\|_{L^{p-1}(\o)} \leq C \|f -
g\|_{L^{q}(\o)}^{\frac 1r} \big(\|f\|_{L^{q}(\o)}+
\|g\|_{L^{q}(\o)}\big)^{\frac 1{p-1} - \frac 1r}
\end{equation*}
for some constant $C$ depending on $p$, $q$ and on the left-hand
side  of either \eqref{+301} or \eqref{+302}. Here, $r = \max \{p,
2\}$.
%
%
%
%
%
%
%
%
%
%
%
%
%
%
%
%
\end{corollary}
Recall from Section \ref{subsec2.3} that inequality \eqref{-514}
between $\lambda$ and $\nu _p$ holds for every domain $\Omega$ and
for every $p>1$, whereas a converse inequality (even up to a
multiplicative constant) fails, unless $\Omega$ is sufficiently
regular. As anticipated in Section \ref{sec1}, Corollaries
\ref{cor3bis} and \ref{cor3} lead to conclusions equivalent to
those of Theorems \ref{teo4bis} and \ref{teo4}, respectively, only
if the domain $\o$ is regular enough for the two sides
%
of \eqref{-514} to be equivalent, namely if a constant $C$ exists
such that
\begin{equation}\label{42}
\nu _p(s) \leq C\bigg(\int _s^{\m2} \frac{dr}{\lambda
(r)^{p'}}\bigg)^{1-p} \qquad \quad \hbox{for $s\in (0, \m2 )$.}
\end{equation}
This is the case of Examples 1--5 below. However, if $\o$ is very
irregular, as in Examples 6 and 7, then \eqref{42} fails, and
Corollaries \ref{cor3bis} and \ref{cor3} are essentially weaker
than Theorems \ref{teo4bis} and \ref{teo4}. \par\noindent In our
examples we shall discuss the problem of existence and uniqueness
of solutions to problem \eqref{1} via Theorem \ref{teo4bis} or
Corollary \ref{cor3bis}; it is implicit that the continuous
dependence on the data will follow under the appropriate strong
monotonicity assumption \eqref{5} by Theorem \ref{teo4} or
Corollary \ref{cor3}, respectively.
\medskip
\par\noindent
{\bf Example 1}. (Lipschitz domains).
\par\noindent
{\rm Assume that $\o $ is a connected and bounded open set with a
Lipschitz boundary, and let $1<p\leq n$. Owing to  \eqref{-516}
and \eqref{-516bis}, condition \eqref{402}  is fulfilled. Thus, by
Theorem \ref{teo4bis}, under assumptions \eqref{2}-\eqref{4}, a
unique approximable solution to problem \eqref{1} exists for any
$f \in L^1(\o)$.
\par\noindent
The same conclusion follows from Corollary \ref{cor3bis}, since
\eqref{42} holds in this case.}
\medskip
\par\noindent
{\bf Example 2}. (H\"older domains).
\par\noindent
Let $\o$ be a connected and bounded open set with a H\"older
boundary with exponent $\alpha \in (0,1)$, and let $1<p<\frac
1\alpha (n-1) +1$. By the Sobolev embedding of \cite{La} and by
the equivalence of \eqref{2001}-\eqref{2002}, we have that
\begin{equation}\label{2010}
\nu _p (s)  \geq C s^{1- \frac {\alpha p}{n-1+\alpha}} \quad
\hbox{for $s \in (0, \mo /2)$,}
\end{equation}
for some positive constant $C$. Owing to Theorem \ref{teo4bis}, a
unique approximable solution to \eqref{1} exists for any $\alpha
\in (0,1)$ and for any $f \in L^1(\o)$.
\par\noindent
On the other hand, by \eqref{2003},
$$\lambda (s) \geq C s^{ \frac {n-1 }{n-1+\alpha}} \quad
\hbox{for $s \in (0, \mo /2)$,}$$ for some positive constant $C$.
Thus, \eqref{42} holds, and  the use of Corollary \ref{cor3bis}
leads to the same conclusion  about solutions to \eqref{1}.
\medskip
\par\noindent
{\bf Example 3}. (John and $\gamma$-John domains).
\par\noindent
Let $\gamma \geq 1$. A bounded open set $\o$ in $\rn$ is called a
$\gamma$-John domain if there exist a constant $c \in (0,1)$ and a
point $x_0 \in \o$ such that for every $x \in \o$ there exists a
rectifiable curve $\varpi : [0, l] \to \o$, parametrized by
arclenght, such that $\varpi (0)=x$, $\varpi (l) = x_0$, and
$${\rm dist}\, (\varpi (r) , \partial \o ) \geq c r^\gamma \qquad
\hbox{for $r \in [0, l]$.}$$ The $\gamma$-John domains generalize
the standard John domains, which correspond to the case when
$\gamma =1$ and arise in connection with the study of holomorphic
dynamical systems and quasiconformal mappings. The notion of John
and $\gamma$-John domain has been used in recent years in the
study of Sobolev inequalities. In particular, a result from
\cite{KM} (complementing \cite{HK}) tells us that if $p \geq 1$
and $1 \leq \gamma \leq \frac p{n-1} +1$, then
$$W^{1,p}(\Omega ) \to L ^{\sigma }(\o),$$
where either $\sigma = \frac{np}{(n-1)\gamma +1-p}$ or $\sigma$ is
any positive number, according to whether $\gamma > \frac
{p-1}{n-1}$ or $\gamma \leq \frac {p-1}{n-1}$. By the equivalence
of \eqref{2001} and \eqref{2002}, one has that
$$\nu _p (s) \geq C s^{\frac {p}{\sigma}} \qquad
\hbox{for $s \in (0, \mo /2)$,}$$
for some positive constant $C$. An application of Theorem
\ref{teo4bis} ensures that a unique approximable solution to
\eqref{1} exists for any $f \in L^q(\o)$ if $q>1$ and $1 \leq
\gamma \leq \frac p{n-1} +1$, and also for $f \in L^1(\o)$
provided that $1 \leq \gamma < \frac p{n-1} +1$.
\par\noindent
It is easily verified, on exploiting \eqref{2003}, that the same
conclusions follow from Corollary \ref{cor3bis} as well.
\medskip
\par\noindent
{\bf Example 4}. (A cusp-shaped domain).
\par
\begin{figure}[htb]

\epsfig{file=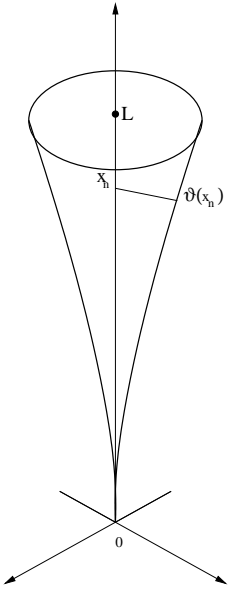, width=3cm} \caption{a  cusp-shaped domain}
\label{cusp}
\end{figure}

\noindent
Let $L>0$ and let  $\vartheta :[0, L] \to [0, \infty )$ be
a differentiable  convex function such that $\vartheta (0) = 0$.
Consider the set
$$ \o = \{x\in \rn : |x'| < \vartheta (x_n ), 0 < x_n <L\}$$
(see Figure \ref{cusp}), where $x=(x', x_n)$ and $x' = (x_1, \dots
, x_{n-1})\in {\R}^{n-1}$. Let $\Theta : [0, L] \to [0, \infty )$
be the function given by
$$\Theta (\rho) = n\omega _{n} \int
_0^\rho \vartheta (r)^{n-1} dr \qquad \quad \hbox{for $\rho \in
[0, L]$.}$$
We claim that \eqref{402} is fulfilled for every $p\in (1,n)$.
Actually, \cite[4.3.5/1]{Ma2} tells us that
\begin{equation}\label{901}
{\nu _p}(s) \approx \bigg(\int _{\Theta
^{-1}(s)}^{\Theta ^{-1}(\m2 )} \vartheta (r)^{\frac{1-n}{p-1}} dr
\bigg)^{1-p} \quad \hbox{for $s \in (0, \m2 )$.}
\end{equation}
Thus, \eqref{402} is equivalent
to
\begin{equation}\label{902}
\int _0s^{-1/p'}\bigg(\int _{\Theta ^{-1}(s)}^{\Theta ^{-1}(\m2 )}
\vartheta (r)^{\frac{1-n}{p-1}} dr \bigg)^{1/p'}ds < \infty\,,
\end{equation}
or, via a change of variable, to
\begin{equation}\label{903}
\int _0\bigg(\int _{\rho}^{\Theta ^{-1}(\m2 )} \vartheta
(r)^{\frac{1-n}{p-1}} dr \bigg)^{1/p'}\bigg(\int _{0}^{\rho}
\vartheta (r)^{n-1} dr\bigg)^{-1/p'}\vartheta (\rho)^{n-1}d\rho <
\infty\,.
\end{equation}
By De L'Hopital rule,
\begin{equation}\label{904}
\limsup _{\rho \to 0^+} \frac{\int _{\rho}^{\Theta ^{-1}(\m2 )}
\vartheta (r)^{\frac{1-n}{p-1}} dr }{\vartheta
(\rho)^{-(n-1)p'}\int _{0}^{\rho} \vartheta (r)^{n-1} dr} \leq
\limsup _{\rho \to 0^+} \frac{1 }{(n-1)p'\vartheta
(\rho)^{-n}\vartheta '(\rho)\int _{0}^{\rho} \vartheta (r)^{n-1}
dr - 1}\,.
\end{equation}
Since $\vartheta '$ is non-decreasing,
\begin{align}\label{905}
\vartheta (\rho)^{-n}\vartheta '(\rho)\int _{0}^{\rho} \vartheta
(r)^{n-1} dr  & = \frac{\vartheta '(\rho) \int _{0}^{\rho}
\big(\int _0^r \vartheta '(\tau) d\tau\big)^{n-1} dr}{\big(\int
_0^\rho \vartheta '(r) dr\big)^{n}}\\
\nonumber & \geq \frac{ \int _{0}^{\rho} \vartheta '(r)\big(\int
_0^r \vartheta '(\tau) d\tau\big)^{n-1} dr}{\big(\int _0^\rho
\vartheta '(r) dr\big)^{n}} =\frac{1}{n}
\end{align}
for $\rho \in (0, \m2)$. Inasmuch as $p<n$, by \eqref{904} and
\eqref{905} the integrand in \eqref{903} is bounded at $0$, and
hence \eqref{903} follows.
\par\noindent
By Theorem \ref{teo4bis}, if $f\in L^1(\o )$, then there exists a
unique approximable solution to problem \eqref{1} under
assumptions \eqref{2}-\eqref{4}.
\par\noindent Notice that the same result can be derived via
Corollary \ref{cor3bis}. Indeed, by \cite[Example 3.3.3/1]{Ma3},
\begin{equation*}
\lambda (s) \approx \vartheta (\Theta ^{-1}(s))^{n-1} \qquad \quad
\hbox{for $s \in (0, \m2 )$,}
\end{equation*}
and hence \eqref{42} holds.
\medskip
\par\noindent
{\bf Example 5} (An unbounded domain).
\par
\begin{figure}[htb]

\epsfig{file=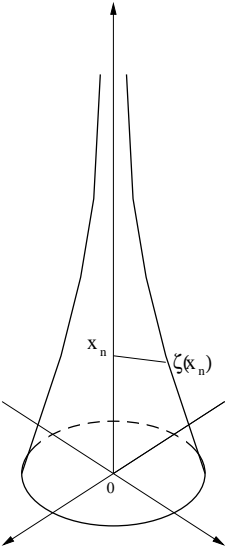, width=3cm} \caption{an unbounded
domain} \label{funnel}
\end{figure}

\noindent
Let $\zeta : [0, \infty )\to (0, \infty )$ be a differentiable
convex function such that $\lim _{\rho \to 0^+} \zeta (\rho) >
-\infty$ and $\lim _{\rho \to \infty } \zeta (\rho) =0$. Consider
the unbounded set
$$ \o = \{x\in \rn : x_n >0, |x'| < \zeta (x_n )\}$$
(see Figure \ref{funnel}), where $x=(x', x_n)$ and $x' = (x_1,
\dots , x_{n-1})\in {\R}^{n-1}$. Assume that
\begin{equation}\label{2020'}
\int_ 0 ^\infty \zeta (r)^{n-1} dr<\infty\,,
\end{equation}
in such a way that $\mo < \infty$. Let $\Upsilon : [0, \infty )
\to [0, \infty )$ be the function given by
$$\Upsilon (\rho) = n\omega _{n} \int
_\rho ^\infty \zeta (r)^{n-1} dr \qquad \quad \hbox{for $\rho
>0$.}$$
By \cite[Example 4.3.5/2]{Ma2}, if $p>1$,
\begin{equation*}
{\nu _p}(s) \approx \bigg(\int _{\Upsilon ^{-1}(\m2 )}^{\Upsilon
^{-1}(s )} \zeta (r)^{\frac{1-n}{p-1}} dr \bigg)^{1-p} \quad
\hbox{for $s \in (0, \m2 )$.}
\end{equation*}
An application of  Theorem \ref{teo4bis} tells us that there
exists a unique solution to problem \eqref{1} with $f \in L^q(\o
)$ if either $q>1$ and
\begin{equation}\label{2020} \int ^\infty
\bigg(\int _{\Upsilon ^{-1}(\m2 )}^\rho \zeta (r)
^{\frac{1-n}{p-1}} dr\bigg)^{\frac {q'}{p'}} \bigg(\int
_\rho^\infty \zeta (r) ^{n-1} dr\bigg)^{\frac {q'}{p}} \zeta
(\rho) ^{n-1}d\rho < \infty\,,
\end{equation}
 or $q=1$ and
\begin{equation}\label{2021}
\int ^\infty \bigg(\int _{\Upsilon ^{-1}(\m2 )}^\rho \zeta (r)
^{\frac{1-n}{p-1}} dr\bigg)^{\frac {1}{p'}} \bigg(\int
_\rho^\infty \zeta (r) ^{n-1} dr\bigg)^{-\frac {1}{p'}} \zeta
(\rho) ^{n-1}d\rho < \infty\,.
\end{equation}
 For
instance, if $\zeta (\rho) = \frac 1{(1+ \rho )^\beta }$, then
\eqref{2020} and \eqref{2020'} hold if $\beta >\frac{1+ q'}{n-1}$,
whereas \eqref{2021} never holds, whatever $\beta $ is. In the
case when $\zeta (\rho ) = e^{-\rho ^\alpha }$ with $\alpha
>0$, condition
 \eqref{2020} holds for every $q \in (1, \infty ]$, whereas \eqref{2021}  does not
 hold for any
$\alpha$.
\par\noindent
Note that, by \cite[Example 3.3.3/2]{Ma2},
$$\lambda (s) \approx \big(\zeta (\Upsilon ^{-1} (s))\big)^{n-1} \qquad
\hbox{as $s\to 0^+$.}$$ Thus, \eqref{2003} holds, and hence
Corollary \ref{cor3bis} leads to the same conclusions.
\medskip
\par\noindent
{\bf Example 6} (A domain from \cite{CH})
\par
\begin{figure}[htb]

\epsfig{file=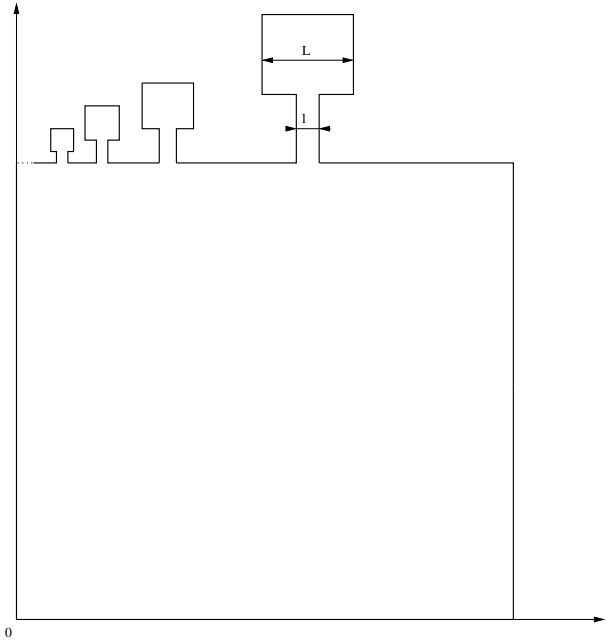, width=8cm} \caption{an example from
\cite{CH}}\label{couhil}
\end{figure}

\noindent
Let us consider problem \eqref{1} in the domain $\o \subset
{\R}^2$ displayed in Figure \ref{couhil} and borrowed from
\cite{CH}, where it is exhibited  as an example of a domain in
which the Poincar\'e inequality fails.
 In the figure, $L = 2^{-k}$ and $l= \delta (2^{- k})$, where $k \in \N$ and $\delta : [0, \infty ) \to [0, \infty )$ is
 any function such that: $\delta (2s) \leq c \delta (s)$ for some $c>0$ and for $s>0$; $\frac{s^{p+1}}{\delta (s)}$ is non-decreasing;
$\frac{s^{1+\varepsilon}}{\delta (s)}$ is non-increasing for some
$\varepsilon >0$. One can show that, if $1\leq p \leq 2$, then
\begin{equation}\label{906}
{\nu _p} (s)  \approx \delta \big( s^{1/2}\big)s^{\frac{1
-p}{2}} \qquad \quad \hbox{as $s \to 0^+$}
\end{equation}
\cite{CM2}. In particular, by \eqref{2003},
\begin{equation}\label{906'}
\lambda  (s)  \approx \delta \big( s^{1/2}\big) \qquad \quad
\hbox{as $s \to 0^+$}.
\end{equation}
By Theorem \ref{teo4bis}, it is easily verified that there exists
a unique solution to problem \eqref{1} if $f\in L^q(\o )$ for any
$q>1$. When $q=1$, the  solution exists and is unique provided
that
\begin{equation}\label{906''}
\int _0\bigg(\frac s{\delta (s)}\bigg)^{1/p} ds < \infty .
\end{equation}
 For instance, \eqref{906''} holds when $\delta (s) = s^\alpha$
 for some $\alpha \in (1, p+1)$, or when $\delta (s) \approx s^{p+1} \big(\log
 (1/s)\big)^{\beta}$ for small $s$, with $\beta >p$.
\par\noindent The use of the isoperimetric function, namely of
Corollary \ref{cor3bis}, yield worses results for the domain of
this example, for which inequality \eqref{42} fails. For instance,
if $\delta (s) = s^\alpha$,
%
%
the existence and uniqueness of a solution to problem \eqref{1}
cannot be deduced from Corollary \ref{cor3bis} unless either
$\alpha < 2$ and $q \geq 1$, or $2 \leq \alpha \leq p+1$ and 
%
%
$q>\frac{2}{4- \delta}$.
\medskip
\par\noindent
{\bf Example 7} (Nikod\'{y}m)
\par
\begin{figure}[htb]
\epsfig{file=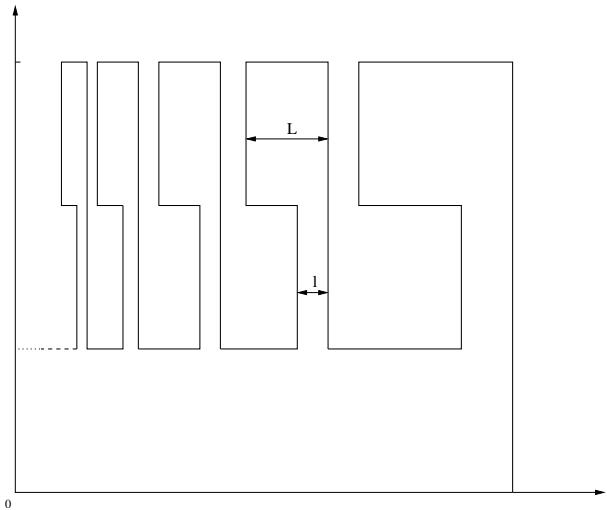, width=8cm} \caption{Nikod\'{y}m
example}\label{niko}
\end{figure}

\noindent
The most irregular domain $\o \subset {\R}^2$ that we consider is
depicted in Figure \ref{niko}. It was introduced by Nikod\'{y}m in
 his study of Sobolev embeddings.
 \par\noindent
 In the figure, $L = 2^{-k}$ and  $l= \delta (2^{- k})$, where $k\in \N$
and
 $\delta : [0, \infty ) \to (0, \infty )$ is
 any increasing Lipschitz continuous function such that $\delta (2s) \leq c \delta (s) \leq c' s$ for some contants $c, c'>0$ and for
 $s>0$.
 If $p\geq 1$, one has that
\begin{equation}\label{909}
{\nu _p} (s) \approx \delta (s)  \qquad \quad \hbox{as $s \to
0^+$},
\end{equation}
and
\begin{equation}\label{908}
\lambda (s) \approx \delta (s)  \qquad \quad \hbox{as $s \to 0^+$}
\end{equation}
\cite[Section 4.5]{Ma2}. By \eqref{909} and Theorem \ref{teo4bis}
there exists a unique approximable solution to problem \eqref{1}
 provided that $f\in L^q(\o )$ for some $q>1$ and \eqref{401} is
 fulfilled, namely
\begin{equation}\label{908'}
\int _0\bigg(\frac{s}{\delta(s)}\bigg)^{\frac {q'}{p}} ds <
\infty\,.
\end{equation}
On the other hand, condition \eqref{402} never holds,
and hence the case when $f\in L^1(\o )$ is not admissible in
Theorem \ref{teo4bis} for this domain.
\par\noindent
In the special case when
\begin{equation}\label{910'}
\delta (s) = s^\alpha, \qquad \quad \hbox{for $s >0$,}
\end{equation}
with $\alpha \geq 1$, condition \eqref{908'} is equivalent to
\begin{equation}\label{910}
\alpha < 1+ \frac {p}{q'} \,.
\end{equation}
\par
Equations \eqref{908} and \eqref{909} tell us that \eqref{42} is
not fulfilled for the domain $\o$ of this example, and the use of
   Corollary
\ref{cor3bis}  actually requires stronger assumtpions on $\delta
(s)$. For instance,  when  $\delta (s)$ is given by \eqref{910'},
one has to demand that $\alpha < 2- \frac 1q$.
%
%
%
This is a more restrictive condition than \eqref{910}.

 {}{}{}{}{}{}
 {}{}{}{}{}{}{}{}{}{}{}{}{}{}{}{}{}{}{}{}{}{}{}{}{}{}{}{}{}{}{}{}{}{}{}{}{}{}
 {}{}{}{}{}{}
 {}{}{}{}{}{}{}{}{}{}{}{}{}{}{}{}{}{}{}{}{}{}{}{}{}{}{}{}{}{}{}{}{}{}{}{}{}{}


\end{document}